\documentclass{article}
\usepackage{amsmath}
\usepackage{amssymb}
\usepackage{mathtools}
\usepackage{verbatim}
\usepackage{enumerate}
\usepackage{graphicx}
\usepackage{float}
\usepackage{overpic}
\usepackage{subfigure}
\usepackage{caption}
\usepackage{geometry}
\usepackage{algorithm}
\usepackage{algorithmic}
\usepackage{multirow}
\usepackage{tikz}
\usepackage{epstopdf}
\usepackage[section]{placeins}
\newcommand\dd{\,\mathrm{d}}
\newcommand\ii{\,\mathrm{i} \,}
\newcommand\bR{\mathbb{R}}

\newcommand\bZ{\mathbb{Z}}

\newcommand\cD{\mathcal{D}}
\newcommand\cI{\mathcal{I}}

\newcommand\pd[2]{\dfrac{\partial{#1}}{\partial{#2}}}
\newcommand{\tdiag}{\text{diag}}

\newcommand{\tif}{\text{if }}

\newcommand{\tfor}{\text{for }}

\begin{document}
\title{Sparsify and sweep: an efficient preconditioner for the Lippmann-Schwinger equation}
\author{Fei Liu$^\sharp$ and Lexing Ying$^{\dagger\sharp}$\\
  $\dagger$ Department of Mathematics, Stanford University\\
  $\sharp$ Institute for Computational and Mathematical Engineering, Stanford University
}
\date{}
\maketitle

\begin{abstract}
This paper presents an efficient preconditioner for the Lippmann-Schwinger equation that combines
the ideas of the sparsifying and the sweeping preconditioners. Following first the idea of the
sparsifying preconditioner, this new preconditioner starts by transforming the dense linear system
of the Lippmann-Schwinger equation into a nearly sparse system. The key novelty is a newly
designed perfectly matched layer (PML) stencil for the boundary degrees of freedoms. The resulting
sparse system gives rise to fairly accurate solutions and hence can be viewed as an accurate discretization of the Helmholtz equation.
This new PML stencil also paves the way for applying the moving PML sweeping preconditioner to
invert the resulting sparse system approximately. When combined with the standard GMRES solver, this
new preconditioner for the Lippmann-Schwinger equation takes only a few iterations to converge for
both 2D and 3D problems, where the iteration numbers are almost independent of the frequency. To the
best of our knowledge, this is the first method that achieves near-linear cost to solve the 3D
Lippmann-Schwinger equation in high frequency cases.
\end{abstract}



\section{Introduction}

\begin{figure}
[!ht]
\centering






\includegraphics[width=\textwidth]{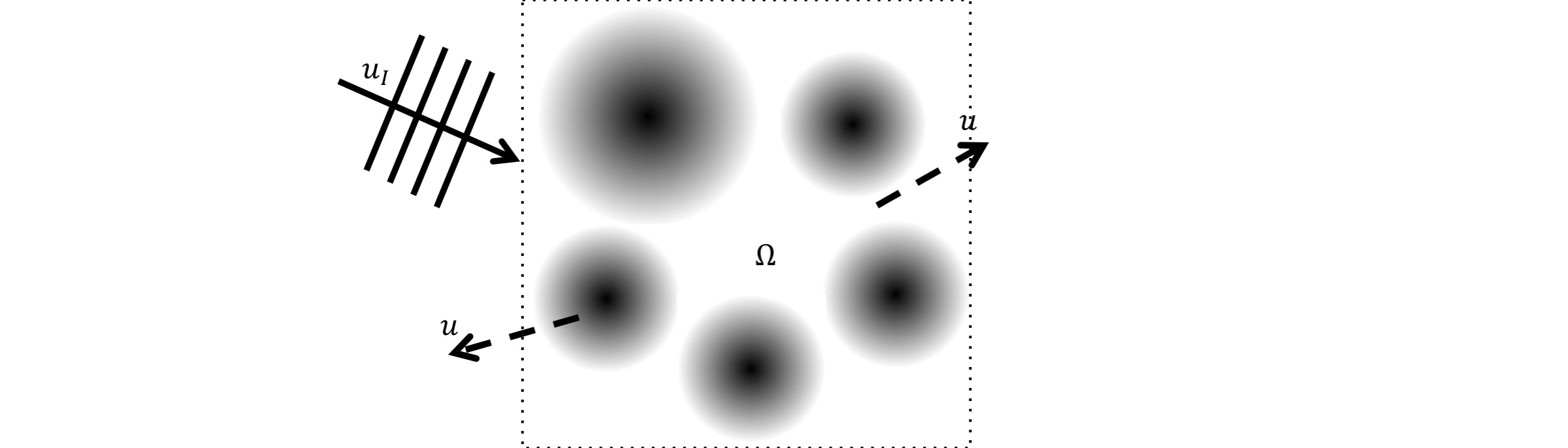}
\caption{An example of the incoming wave $u_I(x)$ and the scattered field $u(x)$.}
\label{fig:scattering}
\end{figure}

This paper concerns the time-harmonic scattering problem
\begin{equation}
\label{eqn:scattered_field}
\begin{dcases}
\left(-\Delta-\dfrac{\omega^2}{c(x)^2}\right)(u(x)+u_I(x)) = 0, \quad x \in \bR^d,\\
\lim_{r\to \infty} r^{(d-1)/2}\left(\pd{}{r}-\ii \omega \right) u(x) = 0,
\end{dcases}
\end{equation}
where $u_I(x)$ is the given incoming wave, $u(x)$ is the scattered field to solve, $\omega$ is the
angular frequency and $c(x)=\Theta(1)$ is the velocity field such that $c(x) = 1$ outside some
bounded region $\Omega$. See Figure \ref{fig:scattering} for an example. The incoming wave $u_I(x)$
satisfies the homogeneous Helmholtz equation
\begin{equation}
  (-\Delta-\omega^2)u_I(x) = 0, \quad x\in \Omega.
\end{equation}

Let $m(x) = 1 - 1/c(x)^2$ be the perturbation field. Rewriting \eqref{eqn:scattered_field} in terms of $m(x)$ we have
\begin{equation}
\label{eqn:scattered_field_m}
(-\Delta-\omega^2 + \omega^2 m(x))u(x) = -\omega^2m(x)u_I(x), \quad x\in \bR^d.
\end{equation}
Let $G(x)$ be the Green's function of the free space Helmholtz equation
\begin{gather*}
G(x)=
\begin{dcases}
\dfrac{\ii}{4}H_0^{(1)}(\omega |x|),\quad & d = 2,\\
\dfrac{\exp(\ii\omega |x|)}{4\pi |x|},\quad & d = 3.
\end{dcases}
\end{gather*}
Convolving both sides of \eqref{eqn:scattered_field_m} with $G(x)$ gives
\begin{equation}
  \label{eqn:LS}
  u(x) + \omega^2 \int_\Omega G(x-y) m(y) u(y) \dd y = -\omega^2 \int_\Omega G(x-y) m(y) u_I(y) \dd y , 
\end{equation}
which is known as the Lippmann-Schwinger equation written in terms of the scattered field $u(x)$.

Solving the integral equation \eqref{eqn:LS} has several advantages compared to solving
\eqref{eqn:scattered_field}. First, since $m(x)$ is compactly supported, we only need to solve
\eqref{eqn:LS} in $\Omega$. The scattered field $u(x)$ for $x\in \Omega^c$ is explicitly given by
\eqref{eqn:LS} once $u(x)$ in $\Omega$ is known. More importantly, the resulting wave field $u(x)$ in
$\bR^d$ automatically satisfies the Sommerfeld radiation condition. On the contrary, for
\eqref{eqn:scattered_field_m} one has to truncate the domain $\bR^d$ to some bounded region and
impose appropriate boundary conditions to simulate the radiation condition. Second, most local
discretizations of \eqref{eqn:scattered_field_m} suffer from the pollution effect \cite{Babuska1997}
due to inaccurate dispersion relations. \eqref{eqn:LS} avoids this problem by leveraging the
Green's function explicitly in the equation.

However, discretizing \eqref{eqn:LS} also raises several issues. First, the resulting linear system
is dense. By the Nyquist theorem, a constant number of points per wavelength is needed to capture
the oscillations, thus the number of degrees of freedom $N$ is at least $\Theta(\omega^d)$. In high
frequency cases, $N$ can be rather large where it is impractical to solve general dense linear
systems with direct method. Second, the discretized system can have very large condition number for non-negligible
perturbations $m(x)$ due to multiple scattering when $\omega$ is large. As a result, most standard iterative solvers require a large number of iterations to converge.

Recently, several progresses have been made to solve the Lippmann-Schwinger equation
\cite{Andersson2005,Bruno2004,Chen2002,greengard2016,vico2016,Lanzara2004,Sifuentes2010,Vainikko2000,Ying2015,Leonardo2016fast}. \cite{vico2016} proposes a numerical scheme that has spectral accuracy for smooth media by truncating the interactions on the physical domain. \cite{greengard2016} presents an adaptive method for the Lippmann-Schwinger equation in 2D.
\cite{Chen2002} solves the 2D Lippmann-Schwinger equation with a technique which is now often
referred to as recursive interpolative factorization or recursive skeletonization, where the setup
cost is $O(N^{3/2})$ and the solve cost is $O(N \log N)$. \cite{Ying2015} approximates the
discretized dense system by a sparse system, and applies the nested dissection factorization
\cite{george1973nested} to the sparse system as a
preconditioner to the original dense system. The costs are dominated by merely the nested dissection
solver, which are $O(N^{3/2})$ and $O(N \log N)$ for setup and solve in 2D, $O(N^2)$ and
$O(N^{4/3})$ for setup and solve in 3D respectively. \cite{Leonardo2016fast} combines the
sparsifying preconditioner \cite{Ying2015} with the method of polarized traces \cite{Demanet2014} to
design a preconditioner for the Lippmann-Schwinger equation in 2D, which achieves $O(N)$ setup and
$O(N \log N)$ solve costs. As far as we know, \cite{Leonardo2016fast} is the first to achieve
near-linear cost in 2D high frequency cases.

Meanwhile, a series of domain decomposition methods were developed to solve the Helmholtz equation
with Sommerfeld radiation condition
\cite{Ying2011a,Ying2011b,Stolk2013,Chen2013a,Chen2013b,Demanet2014,Vion2014,Liu2016}.  The idea is
to divide the domain into slices and impose suitable transmission conditions between these
slices. These methods reduce the computational costs to $O(N)$ for setup and $O(N)$ for solve in 2D,
and $O(N^{4/3})$ for setup and $O(N \log N)$ for solve in 3D, which is a notable improvement over
the nested dissection method. A recursive technique \cite{Liu2015} further reduces both the setup
and solve costs in 3D to $O(N)$.

This work combines the sparsifying preconditioner in \cite{Ying2015} with the sweeping
preconditioner in \cite{Ying2011b} to develop a new preconditioner which solves the
Lippmann-Schwinger equation in near-linear cost. The sketch of the method is as follows. We first
construct two types of compact stencil schemes to approximate the discretized dense system by a
sparse system, and then apply the sweeping factorization to the sparse system. The solving process
of the sweeping factorization induces an approximating solution, which defines a preconditioner to
the original system. The setup and application costs are $O(N)$ and $O(N)$ in 2D and $O(N^{4/3})$
and $O(N \log N)$ in 3D respectively. Furthermore, the costs in 3D can be reduced to $O(N)$ for
setup and $O(N)$ for application by a recursive sweep similar to \cite{Liu2015}. When combined with
the standard GMRES solver, the preconditioner only needs a few iterations to converge, where the
iteration number is almost independent of the angular frequency $\omega$ as shown by the numerical
results. To the best of our knowledge, this is the first algorithm to solve the Lippmann-Schwinger
equation in near-linear cost in 3D high frequency cases.

Another highlight of this work is the newly designed compact stencil introduced for the
preconditioner. The design approach focuses on fitting the stencils to the wave data given by the analytic expressions such as the Green's function. This approach is quite different from the state-of-the-art methods \cite{stolk2016,kristek2009brief,asvadurov2003optimal,babuska1995} to design compact stencils, which focus more on the analytic property of the underlying differential operator. Numerical results show that, when used as a method for solving the Helmholtz equation, this scheme is comparably as accurate as the Quasi-Stabilized FEM (QSFEM) method in \cite{babuska1995} in terms of the phase error.

The rest of the paper is organized as follows. Sections \ref{sec:2D} and \ref{sec:3D} present the
preconditioners and the numerical results in 2D and 3D respectively, where the detailed approach is
explained in Section \ref{sec:2D} for the 2D case, and Section \ref{sec:3D} generalizes it to 3D
with necessary modifications. Section \ref{sec:stencil} presents numerical results to show the
validity of the compact stencil sparsifying scheme presented in this work when used as a direct
method. Conclusions and future work are given in Section \ref{sec:conclusion}.

\section{Preconditioner in 2D}
\label{sec:2D}
This section describes the preconditioner for the 2D Lippmann-Schwinger equation. Starting by
formalizing the dense linear system obtained from discretization, we transform it into an
approximately sparse one by introducing two types of compact stencils. After that, the sweeping
factorization is used to solve the truncated sparse system approximately. The whole process can then
be treated as a preconditioner for the original dense system of the Lippmann-Schwinger equation.

\subsection{Problem formulation}
Without loss of generality, we assume that $\Omega = (0,1)^2$ and that $m(x)$ is supported in
$\Omega$. The task is to discretize the Lippmann-Schwinger equation \eqref{eqn:LS} and solve for
$u(x)$ in $\Omega$.

The domain $\Omega$ is discretized by a uniform Cartesian grid, which allows for the rapid
evaluation of the convolution in \eqref{eqn:LS} by FFT. Let $n$ be the
number of grid points per unit length, $h \coloneqq 1/(n+1)$ be the step size, and $N \coloneqq n^2$
be the number of degrees of freedom.

Denote $i$ as the 2D index point and $p_i$ as the grid point with step size $h$ by
\begin{gather*}
i \coloneqq (i_1, i_2), \quad i_1,i_2 \in \bZ,\\
p_i \coloneqq i h = (i_1 h, i_2 h).
\end{gather*}
Let $\cI$ be the index set of the grid points in $\Omega$ and $\cD$ be the set of the corresponding grid points, given by
\begin{gather*}
\cI \coloneqq \{ i = (i_1,i_2) : 1 \le i_1, i_2 \le n\}, \\
\cD \coloneqq \{ p_i : i \in \cI\}.
\end{gather*}
We also introduce $\bar{\cI}$ as the index set for $\bar{\Omega}$ and $\partial {\cI}$ as the boundary index set by
\begin{gather*}
\bar{\cI} \coloneqq \{i = (i_1,i_2) : 0 \le i_1, i_2 \le n+1\}, \\
\partial \cI \coloneqq \bar{\cI} \setminus \cI,
\end{gather*}
and correspondingly we have $\bar{\cD}$ and $\partial\cD$ as
\begin{gather*}
\bar{\cD} \coloneqq \{ p_i : i \in \bar \cI\},\\
\partial{\cD} \coloneqq \{ p_i : i \in \partial \cI\}.
\end{gather*}

Let $u_i$ be the numerical solution of \eqref{eqn:LS} at $p_i$ for $i\in \cI$. To compute the integral in \eqref{eqn:LS}, we use the Nystr\"om method
\begin{gather*}
\int_\Omega G(p_i-y) m(y) u(y) \dd y \approx \sum_{j \in \cI} k_{i - j} m_j u_j,
\end{gather*}
where
\begin{gather*}
m_i \coloneqq m(p_i), \quad i \in \cI,\\
k_i \coloneqq G(p_i)h^2, \quad i \ne (0,0),
\end{gather*}
and $k_{(0,0)}$ is the weight given by a quadrature correction at the singular point of $G(x)$ at $x=0$, which achieves $O(h^4 \log(1/h)^2)$ accuracy when $m(x)$ is smooth \cite{Duan2009}. This gives the discretized equation
\begin{equation}
\label{eqn:LS_d}
u_i + \omega ^2 \sum_{j \in \cI } k_{i - j} m_j u_j = g_i, \quad  i \in \cI,
\end{equation}
where
\begin{gather*}
g_i \coloneqq - \omega^2 \sum_{j \in \cI} k_{i - j} m_j [u_I]_j, \quad i\in \cI,
\end{gather*}
and $[u_I]_j \coloneqq u_I(p_j)$ is the discrete value of the incoming wave. Higher order quadrature
can be achieved by using more extended local quadrature correction \cite{Duan2009}.

With a slight abuse of the notations, we extend the discrete vectors $m$ and $g$ to the whole 2D grid by zero padding
\begin{gather*}
m_i \coloneqq 0,\quad i \in \bZ^2 \setminus \cI, \\
g_i \coloneqq 0,\quad i \in \bZ^2 \setminus \cI.
\end{gather*}
Introducing matrix $K$ with $K_{i,j} \coloneqq k_{i-j}$, \eqref{eqn:LS_d} can be written into a more compact form
\begin{equation}
\label{eqn:LS_d_c}
(I + \omega^2 K M) u = g,
\end{equation}
where $M \coloneqq \tdiag(m)$.

A subtle difference between \eqref{eqn:LS_d} and \eqref{eqn:LS_d_c} is that, \eqref{eqn:LS_d} is a
set of equations for the unknowns with indices $i \in \cI$, while \eqref{eqn:LS_d_c} can be regarded
as an equation set defined on the infinite index set $\bZ^2$, where the unknown vector $u$ is also
extended to the whole 2D grid with the extension value determined by the equation implicitly. We
have two observations for \eqref{eqn:LS_d_c}
\begin{itemize}
\item
  The solution of \eqref{eqn:LS_d_c} agrees with the one of \eqref{eqn:LS_d} in $\cI$. To get the
  numerical solution of \eqref{eqn:LS} in $\Omega$, we can solve \eqref{eqn:LS_d_c} and then
  restrict the solution to $\cI$ instead of solving \eqref{eqn:LS_d}.
\item
  The solution of \eqref{eqn:LS_d_c} does not match the numerical solution of \eqref{eqn:LS} outside
  $\Omega$ since the zero padding of $g$ differs from the discretized value of the right-hand side of
  \eqref{eqn:LS} in $\Omega^c$. Nonetheless, this is not an issue as we only care about the solution
of \eqref{eqn:LS} in $\Omega$.
\end{itemize}

One may ask: why do we extend the discrete domain to the infinite grid and consider a problem with
infinite size? Besides, the zero padding pattern of $g$ seems rather irrational as it creates
discontinuities at $\partial \cD$. The answer is that, we are not going to actually solve the
$\bZ^2$-size problem. The purpose of extending the unknown to a larger domain is to introduce the
wave attenuation by PML on the extended grid to simulate the Sommerfeld radiation condition as we
shall see in Section \ref{subsubsec:PML_stencil}. The zero padding of $g$ is to ensure that there is
no source outside $\Omega$ such that the PML approximation holds.

The reader may notice that, if we just use the discretized value of the right-hand side of
\eqref{eqn:LS} defined on the whole plane, the solution will also satisfy the Sommerfeld condition, so
it seems meaningless to introduce the zero padding. It is true that the right-hand side of
\eqref{eqn:LS} on the whole plane will induce a solution satisfying the radiation
condition. However, in some cases, when solving \eqref{eqn:LS_d}, we are only given $g$ defined in
$\cI$ without knowing the actual incoming wave $u_I$, and it's computationally impractical to get
the extension of $g$ determined by \eqref{eqn:LS}. This is especially true when we develop
preconditioners where the input only involves the right-hand side in the domain of interest.

With the extended problem \eqref{eqn:LS_d_c}, we are now ready to build a sparse system to
approximate \eqref{eqn:LS_d}.

\subsection{Sparsification}
\label{subsec:2D_sp}
In this section, we adopt the idea of the sparsifying preconditioner \cite{Ying2015} to build a
sparse system which serves as an approximation to \eqref{eqn:LS_d}. The sparse system to be
constructed has the same sparsity pattern as a compact stencil scheme, i.e., each equation only
involves the unknowns at one grid point and its neighbor points, unlike \eqref{eqn:LS_d} where each
equation is dense in $\cI$.

To be specific, we define $\mu_i$ as the neighborhood for the index $i$
\begin{gather*}
\mu_i \coloneqq \{j : \|j - i\|_\infty \le 1\}.
\end{gather*}
Now the task is to build for each point $i$ a local stencil supported only in $\mu_i$.  We shall
build two types of stencils in what follows. The first type is for the interior points, while the second type is for the points near the boundary which are inside
what we call ``the PML region''.

The perfectly matched layer (PML) \cite{Berenger1994,Johnson2008,Chew1994} is a technique to
attenuate the waves exponentially near the boundary of the domain so that the zero Dirichlet
boundary conditions can be imposed directly to simulate the radiation boundary condition without
bringing in too much error. We will explain the PML usage during the construction of the second type
of the stencils.

Let's start with extending our domain $\Omega$. Denote $\Omega^{h}$ as $\Omega$ with an $h$-size extension, $\Omega^{h+\eta}$ as the PML extension of $\Omega^h$ with width $\eta$, given by
\begin{gather*}
\Omega^h \coloneqq (-h, 1+h)^2,\\
\Omega^{h+\eta} \coloneqq (-h -\eta, 1+h + \eta)^2.
\end{gather*}
Here $\eta = bh$ is the PML width, where $b = O(1)$ is the number of discrete layers in each side of
the PML region. $\eta$ is typically around one wavelength. The PML region $\Omega^{h+\eta} \setminus
\Omega^{h}$ is where we will attenuate the scattered field $u(x)$ in Section
\ref{subsubsec:PML_stencil}. Note that there is a small $h$-distance between the domain of interest
$\Omega$ and the PML region $\Omega^{h+\eta} \setminus \Omega^{h}$. This small distance is
introduced on purpose and the reason will be clear later.

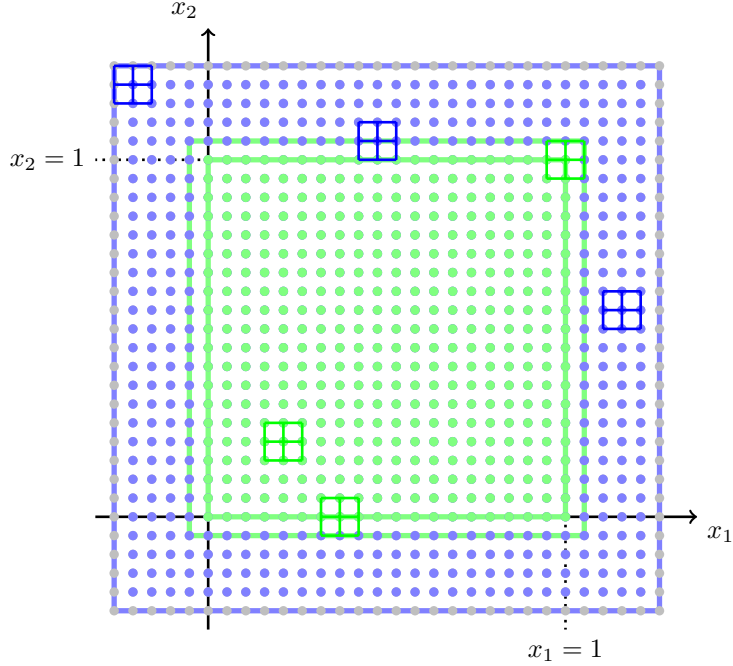
\begin{figure}
[!ht]
\centering

\begin{tikzpicture}
[scale = 0.25]
\draw[->,line width = 1pt](-6,0) -- (26,0)node[below right]{$x_1$};
\draw[->,line width = 1pt](0,-6) -- (0,26)node[above left]{$x_2$};

\draw[dotted, line width = 1pt](19,0)--(19,-6) node [below] {$x_1=1$};
\draw[dotted, line width = 1pt](0,19)--(-6,19) node [left] {$x_2=1$};

\draw(0,0)[line width = 2pt,green!50]rectangle(19,19);
\draw(-1,-1)[line width = 2pt,green!50]rectangle(20,20);
\draw(-5,-5)[line width = 2pt,blue!50]rectangle(24,24);

\foreach \x in {-5,...,24}
\foreach \y in {-5,...,24}
\fill(\x,\y)[fill=gray!50]circle(0.25);

\foreach \x in {-4,...,23}
\foreach \y in {-4,...,23}
\fill(\x,\y)[fill=blue!50]circle(0.25);

\foreach \x in {0,...,19}
\foreach \y in {0,...,19}
\fill(\x,\y)[fill=green!50]circle(0.25);

\draw (3,3)[green!100, line width = 1pt]grid(5,5);
\draw (6,-1)[green!100, line width = 1pt]grid(8,1);
\draw (18,18)[green!100, line width = 1pt]grid(20,20);

\draw (21,10)[blue!100, line width = 1pt]grid(23,12);
\draw(8,19)[blue!100, line width = 1pt]grid(10,21);
\draw(-5,22)[blue!100, line width = 1pt]grid(-3,24);

\end{tikzpicture}
\caption{This figure is an illustration of $\Omega, \Omega^h, \Omega^{h+\eta}, \cD, \cD^{h}$, etc where $n=18, b=4$. The inner green box is the boundary of $\Omega$, the outer green box (with blue points on it) is the boundary of $\Omega^h$ and the big blue box (with gray points on it) is the boundary of $\Omega^{h+\eta}$. The set of the all the green points is $\cD^h$, and the set of all the blue points is $\cD^{h+\eta}\setminus \cD^h$. $\cD$ is the set of all the green points strictly inside the inner green box, while $\partial \cD$ consists of the points exactly located on the inner green box. The gray points form $\partial \cD^{h+\eta}$ where we will impose the zero Dirichlet boundary conditions after wave attenuation. Each of the $3\times 3$ green grid corresponds to a neighborhood $\mu_i$ for some $i \in \cI^h$ where we construct stencils of the first type. Correspondingly, the $3\times 3$ blue grids are for $\mu_i$ with $i \in \cI^{h+\eta} \setminus \cI^h$ where we build stencils of the second type.}
\label{fig:2D_d}
\end{figure}

The corresponding index sets in these regions are
\begin{gather*}
\cI^h \coloneqq \{i : 0 \le i_1, i_2 \le n + 1\},\\
\cI^{h+\eta} \coloneqq \{ i : -b \le i_1, i_2 \le n + 1 + b\}.
\end{gather*}
Similar to the notations $\bar\cI, \partial \cI, \cD, \partial \cD$ and $\bar{\cD}$, we introduce
$\bar\cI^h, \partial \cI^h, \cD^h, \cD^{h+\eta}, \partial \cD^h$, etc, as the corresponding grid
point sets, boundary sets, closures and so on. The meanings are straightforward and we omit the
formal definitions. See Figure \ref{fig:2D_d} for an illustration.

We now describe how to design two types of stencils for the unknowns indexed by $\cI^{h+\eta}$: first
for the ones in $\cI^h$ and then for the ones in $\cI^{h+\eta}\setminus \cI^h$. At the end, we
assemble them together to form our sparse system.

\subsubsection{Stencils for the interior points in $\cI^h$}
\label{sec:Ih_stencil}
Following the approach in \cite{Ying2015}, we design the first type of the stencils for the
neighborhood $\mu_i$ where $i \in \cI^h$ (see Figure \ref{fig:2D_d}, the $3\times 3$ green
grids). Taking out the equations in \eqref{eqn:LS_d_c} indexed by $\mu_i$ we have
\begin{equation}
u_i + \omega ^2 \sum_{j \in \cI } K_{i,j} m_j u_j = g_i, \quad  i \in \mu_i,
\end{equation}
which can be written as
\begin{equation}
\label{eqn:local}
u_{\mu_i} + \omega^2 (K_{\mu_i,\mu_i} [m u]_{\mu_i} + K_{\mu_i,\mu_i^c} [mu]_{\mu_i^c}) = g_{\mu_i}.
\end{equation}
Here are some explanations for the notations in \eqref{eqn:local}:
\begin{itemize}
\item
  The subscript $\mu_i$ stands for the corresponding vector restricted to the index set $\mu_i$, for
  example, $[mu]_{\mu_i}$ is the vector of the elementwise multiplication of $m$ and $u$ restricted
  to $\mu_i$.
\item
  $\mu_i^c \coloneqq \cI \setminus \mu_i$, which is the complement of $\mu_i$ with respect to $\cI$.
\item
  $K_{\mu_i, \mu_i^c}$ is the sub-matrix of $K$ with row index set $\mu_i$ and column index set $\mu_i^c$.
\end{itemize}

Let's consider a linear combination of the equations in \eqref{eqn:local}. Suppose $\alpha$ is a
column vector supported on $\mu_i$. Multiplying both sides of \eqref{eqn:local} by $\alpha^*$ gives
\begin{equation}
\label{eqn:local_lc}
\alpha^* u_{\mu_i} + \omega^2 (\alpha^* K_{\mu_i,\mu_i} [m u]_{\mu_i} + \alpha^* K_{\mu_i,\mu_i^c} [mu]_{\mu_i^c}) = \alpha^* g_{\mu_i},
\end{equation}
where $\alpha^*$ is the conjugate transpose of $\alpha$.

To design a local stencil, we hope that the resulting equation \eqref{eqn:local_lc} only involves unknowns indexed by $\mu_i$. Observing the left-hand side of \eqref{eqn:local_lc}, we found if $\alpha^* K_{\mu_i,\mu_i^c} \approx 0$, then we can truncate the terms involving $u_{\mu_i^c}$ and the resulting equation will be local. But does there exist an $\alpha$ such that $\alpha^* K_{\mu_i,\mu_i^c} \approx 0$? The answer is yes. The reason is that the elements of $K$ are defined by the Green's function $G(x)$, which satisfies
\begin{equation}
(-\Delta - \omega^2) G(x) = 0,\quad  x \in \bR^2 \setminus \{(0,0)\}.
\end{equation}
Each column of the matrix $K_{\mu_i,\mu_i^c}$ can be treated as the Green's function centered at some grid point indexed by $j \in \mu_i^c$ and evaluated at the points indexed by the neighborhood $\mu_i$, which does not involve the singular point of $G(x)$ at $x= 0$. Thus it's reasonable to expect some local stencil $\alpha$, which can be thought of as a discretization of the local operator $(-\Delta - \omega^2)$, such that $\alpha^* K_{\mu_i,\mu_i^c} \approx 0$. By the translational invariance of the Green's function, to find such $\alpha$, it suffices to require that $\alpha^* K_{\mu, \mu^c} \approx 0$, where
\begin{gather*}
\mu \coloneqq \mu_0 = \{j : \|j\|_\infty \le 1 \},\\
\mu^c \coloneqq \{i : -n \le i_1, i_2 \le n\} \setminus \mu,
\end{gather*}
which means that we can translate the index $i$ to the origin and consider an equivalent problem. Here the complement of $\mu$ is taken with respect to a larger index set. The reason is that, when we translate different indices $i$ to the origin, the corresponding complement $\mu_i^c$ will also be translated. The larger set is taken as the union of all those translated complements to ensure that the condition is sufficient.

To minimize $\alpha^* K_{\mu, \mu^c}$, we consider the optimization problem
\begin{equation}
  \label{eqn:alpha}
  \min_{\alpha : \|\alpha\|_2} \|\alpha^* K_{\mu, \mu^c}\|_2.
\end{equation}
The solution is the left singular vector corresponding to the smallest singular value of $K_{\mu, \mu^c}$, which can be solved in $O(N)$.

Once we have $\alpha$, we compute $\beta$ by setting
\begin{equation}
  \label{eqn:beta}
  \beta^* \coloneqq \alpha^* K_{\mu, \mu}.
\end{equation}
Then \eqref{eqn:local_lc} can be approximated as
\begin{equation}
  \label{eqn:local_c_approx}
\alpha^* u_{\mu_i }+ \omega^2 \beta^* [m u]_{\mu_i} \approx \alpha^* g_{\mu_i}.
\end{equation}
This defines the local stencil for each $i \in \cI^h$.

Note that, if we do the same thing for $i \notin \cI^h$, the right-hand side $\alpha^* g_{\mu_i}$ will be $0$ due to the zero padding of $g$. If we build the stencils for all $i \in \bZ^2 \setminus \cI^h$ and combine them with the Sommerfeld radiation condition at infinity, it will induce a discrete DtN map at $\partial \cI^h$. This linear map, though existing in theory, is dense and expensive to compute. Section \ref{subsubsec:PML_stencil} circumvents this issue by exploiting PML on the extended domain and introducing the second type of stencils to approximate this dense map efficiently.

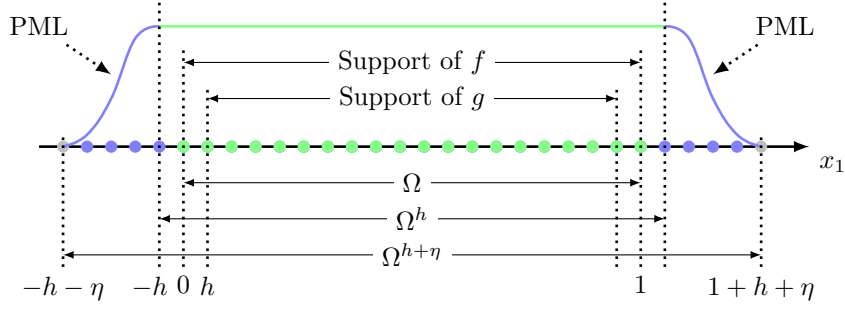
\begin{figure}
[!ht]
\centering

\begin{tikzpicture}
[>=latex,scale=0.32]

\draw[->,line width = 1pt](-6,0) -- (26,0)node[below right]{$x_1$};

\foreach \x in {-5,...,24} {
\fill (\x,0) [fill=gray!50]circle(0.25);
}

\foreach \x in {-4,...,23} {
\fill (\x,0) [fill=blue!50]circle(0.25);
}

\foreach \x in {0,...,19} {
\fill (\x,0) [fill=green!50]circle(0.25);
}

\draw [green!50, line width = 1pt] (-1,5) -- (20,5);
\draw [blue!50, line width = 1pt] plot [smooth, tension = 0.7] coordinates {(20, 5) (21,4.5) (22,2) (23,0.5) (24,0)};
\draw [blue!50, line width = 1pt] plot [smooth, tension = 0.7] coordinates {(-1, 5) (-2,4.5) (-3,2) (-4,0.5) (-5,0)};

\draw [dotted, line width = 1pt] (1, -5) -- (1, 2.5);
\draw [dotted, line width = 1pt] (18, -5) -- (18, 2.5);

\draw [dotted, line width = 1pt] (0, -5) -- (0, 4);
\draw [dotted, line width = 1pt] (19, -5) -- (19, 4);

\draw [dotted, line width = 1pt] (-1, -5) -- (-1, 6);
\draw [dotted, line width = 1pt] (20, -5) -- (20, 6);

\draw [dotted, line width = 1pt] (-5, -5) -- (-5, 0.5);
\draw [dotted, line width = 1pt] (24, -5) -- (24, 0.5);

\draw[<->] (1, 2) -- (18, 2) node [midway, fill=white] {Support of $g$};
\draw[<->] (0, 3.5) -- (19, 3.5) node [midway, fill=white] {Support of $f$};

\draw[<->] (0, -1.5) -- (19, -1.5) node [midway, fill=white] {$\Omega$};
\draw[<->] (-1, -3) -- (20, -3) node [midway, fill=white] {$\Omega^h$};
\draw[<->] (-5, -4.5) -- (24, -4.5) node [midway, fill=white] {$\Omega^{h+\eta}$};

\draw (1, -5) node [below] {$h$};
\draw (0, -5) node [below] {$0$};
\draw (-1.4, -5) node [below] {$-h$};
\draw (-5, -5) node [below] {$-h-\eta$};

\draw (19, -5) node [below] {$1$};
\draw (24, -5) node [below] {$1+h+\eta$};

\draw (25, 5) node (PML right) {PML};
\draw [->, dotted, line width = 1pt] (PML right) -- (22, 3);

\draw (-6, 5) node (PML left) {PML};
\draw [->, dotted, line width = 1pt] (PML left) -- (-3, 3);

\end{tikzpicture}

\caption{This figure is an illustration of the supports of the discrete vectors $g$ and $f$ in 1D where $n= 18, b =4$ and  $f_i \coloneqq \alpha^* g_{\mu_i}$. We see that by introducing the stencil $\alpha$, the support of $g$ is enlarged to the support of $f$ by one grid point on each side.}
\label{fig:1D_d}
\end{figure}

Now why do we introduce the $h$-size padded domain $\Omega^h$ and build the first type of stencils for $i \in \cI^h$ rather than just for $\cI$? The reason is that, for $i\in \partial \cI$, $\alpha^* g_{\mu_i}$ is not necessarily zero, thus we cannot assign $\partial \cI$ to the second type where the corresponding right-hand side is zero. So we enlarge $\Omega$ by $h$-size and build stencils of the first type for $\cI^h = \bar{\cI}$. Figure \ref{fig:1D_d} shows this subtlety in 1D as an illustration.

\subsubsection{Stencils for the PML points in $\cI^{h+\eta}\setminus \cI^h$}
\label{subsubsec:PML_stencil}
Next, we design the stencils for $i \in \cI^{h+\eta}\setminus \cI^h$ (see Figure \ref{fig:2D_d}, the $3\times 3$ blue grids). Define the auxiliary function
\begin{gather*}
\sigma(x) \coloneqq
\begin{dcases}
- \dfrac{C}{\omega} \left(\dfrac{x + h}{\eta}\right)^2, \quad & \tif -h - \eta <  x \le -h,\\
0, \quad & \tif -h < x < 1 + h,\\
\dfrac{C}{\omega} \left(\dfrac{x - 1 -h }{\eta}\right)^2, \quad &  \tif 1+h \le x < 1 + h + \eta,
\end{dcases}
\end{gather*}
where $C \sim \Theta(1)$ is some positive constant. We attenuate the scattered field $u(x)$ in the PML region $\Omega^{h+\eta}\setminus \Omega^h$ by introducing the complex stretching
\begin{gather*}
x^\sigma \coloneqq (x^\sigma_1,x^\sigma_2) = (x_1 + \ii \sigma(x_1), x_2 + \ii \sigma(x_2)),\\
u^\sigma (x) \coloneqq u(x^\sigma) = u(x_1 + \ii \sigma(x_1), x_2 + \ii \sigma(x_2)),\\
u^\sigma_i \coloneqq u^\sigma (p_i) = u(p^\sigma_i).
\end{gather*}

By changing variable from $x$ to $x^\sigma$, we know that the function $u^\sigma(x)$ satisfies the modified Helmholtz equation in the PML region

\begin{equation}
\label{eqn:H_PML}
\left(-\sum_{d = 1} ^ 2 \left(\dfrac{\partial_d}{1+ \ii \sigma'(x_d)}\right)^2 - \omega^2 \right) u^\sigma (x) = 0, \quad  x \in \Omega^{h+\eta}\setminus \Omega^h,
\end{equation}

A simple way to build local stencils for $\cI^{h+\eta} \setminus \cI^h$ is to discretize
\eqref{eqn:H_PML} explicitly with some local scheme such as the central difference
scheme. Unfortunately, it turns out to be not accurate enough to do so. We adopt a different
approach. The idea is similar to what we did in the previous section, where we aim to find some
local stencil to annihilate a set of given functions evaluated at the points indexed by $\mu_i$. In
what above, we used the Greens function $G(x)$ to design the stencil $\alpha$. Here we use a set of
``modified plane waves'' to achieve the same goal.

Specifically, we first note that the plane wave function
\begin{gather*}
  F(x) \coloneqq \exp(\ii \omega (r\cdot x)), \quad \|r\|_2 = 1,
\end{gather*}
satisfies the free space Helmholtz equation
\begin{equation}
  \label{eqn:H}
  \left(-\Delta - \omega^2 \right) F(x) = 0, \quad  x \in \bR^2.
\end{equation}
Let $F^\sigma(x) \coloneqq F(x^\sigma)$ be the complex stretching of $F(x)$. We immediately have
that $F^\sigma(x)$ satisfies \eqref{eqn:H_PML} by definition. If we were to design a local stencil
$\gamma $ for $\mu_i$ where $i\in \cI^{h+\eta} \setminus \cI^h$, we would hope that $\gamma^*
F^\sigma_{\mu_i} \approx 0$, where $F^\sigma_{\mu_i}$ is the function $F^\sigma(x)$ evaluated at the
grid points indexed by $\mu_i$. Note that any direction $r$ such that $\|r\|_2=1$ induces a
``modified plane wave'' $F^\sigma(x)$. We hope to solve $\gamma$ by annihilating as many $r$ as
possible. To be precise, Let $R$ be a set of directions where the elements are sampled uniformly
from the unit circle $\{r: \|r\|_2=1\}$, and $F^\sigma_{\mu_i, R}$ be a matrix of size $|\mu_i|
\times |R|$, each column of which is a modified plane wave function $F^\sigma(x)$ with a direction
$r\in R$, evaluated at the grid points indexed by $\mu_i$. Then we solve $\gamma$ by
\begin{equation}
  \label{eqn:gamma}
  \min_{\gamma:\|\gamma\|_2 = 1} \| \gamma^* F^\sigma_{\mu_i, R} \|_2.
\end{equation}

Intuitively, it's better to increase the sample size $|R|$ to improve the reliability of the
stencil. However, larger sample size also leads to more computational cost. Fortunately, it turns
out that not too many samples are needed for a reliable result. It suffices to use only the eight
most common directions -- north, south, west, east, northwest, northeast, southwest and southeast --
to form $R$, and $\gamma$ is given by the vector perpendicular to the eight corresponding vectors on
$\mu_i$. Note that the solution is unique up to a coefficient $\pm 1$ since we have $8$ independent
modified plain waves and the size of the neighborhood $\mu_i$ is $9$.

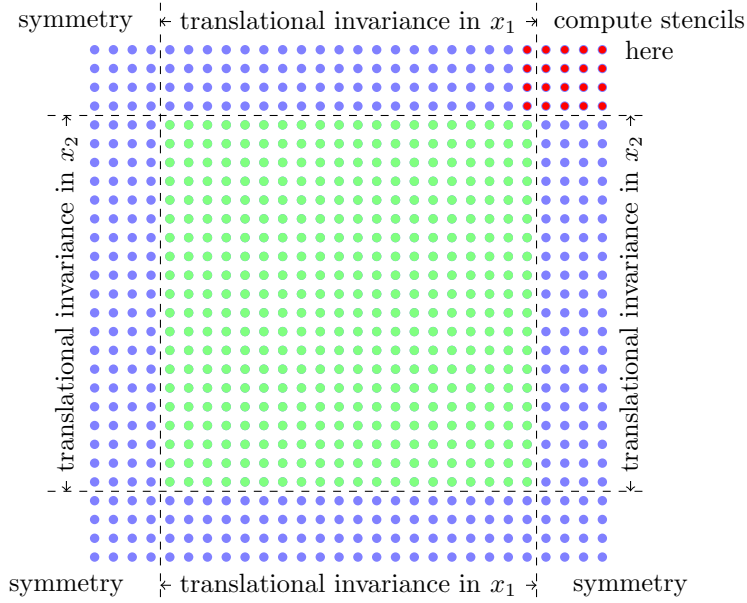
\begin{figure}
[!ht]
\centering

\begin{tikzpicture}
[scale = 0.25]

\foreach \x in {-4,...,23}
\foreach \y in {-4,...,23}
\fill(\x,\y)[fill=blue!50]circle(0.25);

\foreach \x in {0,...,19}
\foreach \y in {0,...,19}
\fill(\x,\y)[fill=green!50]circle(0.25);

\draw[dashed](-0.5,25.5)--(-0.5,-6.5);
\draw[dashed](19.5,25.5)--(19.5,-6.5);
\draw[dashed](-6.5,-0.5)--(25.5,-0.5);
\draw[dashed](-6.5,19.5)--(25.5,19.5);

\draw[<->](-0.5,-5.5)--(19.5,-5.5)node[midway,fill=white]{translational invariance in $x_1$};
\draw[<->](-0.5,24.5)--(19.5,24.5)node[midway,fill=white]{translational invariance in $x_1$};

\draw[<->](-5.5,-0.5)--(-5.5,19.5)node[midway,fill=white,rotate=90]{translational invariance in $x_2$};
\draw[<->](24.5,-0.5)--(24.5,19.5)node[midway,fill=white,rotate=90]{translational invariance in $x_2$};

\draw(-5.5,-5.5)node{symmetry};
\draw(-5,24.5)node{symmetry};
\draw(24.5,-5.5)node{symmetry};
\draw(25.5,24.5)node{compute stencils}(25.5,23)node{here};

\foreach \x in {19,...,23}
\foreach \y in {20,...,23}
\fill(\x,\y)[fill=red!100]circle(0.20);

\end{tikzpicture}
\caption{By the symmetry and the translational invariance of the complex stretching along each dimension, one only needs to compute the PML stencils for the points marked with red color near the top right corner.}
\label{fig:2D_PML}
\end{figure}

In the PML region, we need to compute different stencils for different neighborhoods due to the lack of translational invariance as a result of the complex stretching. Nevertheless, by the symmetry of the stretching, we only need to compute the stencils near a corner of $\cI^{h+\eta} \setminus \cI^h$, which takes only $O(b^2)$ work in total. See Figure \ref{fig:2D_PML} for an illustration.

We denote $\gamma_i$ as the stencil for $\mu_i$, then the corresponding approximating equation is
\begin{equation}
\label{eqn:PML_stencil}
\gamma_i^* u^\sigma_{\mu_i} \approx 0.
\end{equation}
This defines the local stencil for each $i \in \cI^{h+\eta} \setminus \cI^h$.

\subsubsection{Assemble together}
Assembling \eqref{eqn:local_c_approx} and \eqref{eqn:PML_stencil} together and noting that
$u^\sigma_i = u_i$ for $i\in \bar{\cI}^h$, we have
\begin{equation}
  \label{eqn:sp_stc}
  \begin{dcases}
    \alpha^* u^\sigma_{\mu_i }+ \omega^2 \beta^* [m u^\sigma]_{\mu_i} \approx \alpha^* g_{\mu_i},\quad & i \in \cI^h,\\
    \gamma_i^* u^\sigma_{\mu_i} \approx 0, \quad & i \in \cI^{h+\eta} \setminus \cI^h,
  \end{dcases}
\end{equation}
where $\alpha$, $\beta$, and $\gamma$ are given in \eqref{eqn:alpha}, \eqref{eqn:beta} and
\eqref{eqn:gamma} respectively. Noticing also that $u^\sigma$ almost satisfies the zero Dirichlet
boundary conditions
\begin{gather*}
  u^\sigma_i \approx 0, \quad i \in \partial \cI^{h+\eta},
\end{gather*}
we can introduce the sparse linear system
\begin{gather*}
  \begin{dcases}
    \alpha^* \tilde u_{\mu_i }+ \omega^2 \beta^* [m \tilde u]_{\mu_i} = \alpha^* g_{\mu_i},\quad & i \in \cI^h,\\
    \gamma_i^* \tilde u_{\mu_i} = 0, \quad & i \in \cI^{h+\eta} \setminus \cI^h,\\
    \tilde u_i = 0, \quad & i \in \partial \cI^{h+\eta}.
  \end{dcases}
\end{gather*}
for the unknown $\tilde u$ defined on $\cD^{h+\eta}$ that serves as an approximation to $u^\sigma$.
In what follows, we write this system conveniently as 
\begin{equation}
  \label{eqn:Huf}
  H \tilde{u} = f,
\end{equation}
where the right-hand side $f$ is given by
\begin{gather*}
  f_i \coloneqq
  \begin{dcases}
    \alpha^* g_{\mu_i},\quad &\tif i \in \cI^h,\\
    0,\quad & \tif i \in \cI^{h+\eta} \setminus \cI^h.
  \end{dcases}
\end{gather*}

The system \eqref{eqn:Huf} is defined on $\cD^{h+\eta}$. To get the unknowns on $\cD$, we simply
solve \eqref{eqn:Huf} and extract the solutions on $\cD$. The result is an approximation to the true
solution of \eqref{eqn:LS_d}, and this process can serve as a preconditioner for the linear system
\eqref{eqn:LS_d}. In the next section, we present an approach for approximating the solution of \eqref{eqn:Huf} efficiently by leveraging the idea of the sweeping preconditioner.

\subsection{Sweeping factorization}
In this section, we adopt the sweeping factorization to solve the sparse system \eqref{eqn:Huf}
approximately. The main idea of the sweeping factorization is to divide the domain into slices and
eliminate the unknowns slice by slice. An auxiliary PML region is introduced for each slice to build
a subproblem to approximate the inverse of the Schur complement during the Gaussian elimination to
save computational cost.

To be specific, we first divide the 2D grid into $\ell$ slices along the $x_1$ direction. Each slice
contains only a few layers. The leftmost slice contains the left PML region and the rightmost one
contains the right PML region (see Figure \ref{fig:2D_slice}). For simplicity, we assume that each
of the middle slices contains $b$ layers and each of the two boundary slices contains $2b$ layers --
$b$ normal layers plus $b$ attenuating layers in the PML region. Let $\cD_1,\dots, \cD_\ell$ be the
discrete points in each slice correspondingly, and define $\tilde u_{[i]}$ and $f_{[i]}$ as the
restrictions of $\tilde u$ and $f$ on $\cD_i$ respectively. The sparse system \eqref{eqn:Huf} can be
written as the block tridiagonal form
\begin{gather*}
\begin{bmatrix}
H_{[1,1]} & H_{[1,2]} \\
H_{[2,1]} & H_{[2,2]} & \ddots \\
 & \ddots & \ddots & \ddots \\
 & & \ddots & \ddots & H_{[\ell-1, \ell]} \\
 & & & H_{[\ell, \ell-1]} & H_{[\ell, \ell]}
\end{bmatrix}
\begin{bmatrix}
\tilde u_{[1]} \\
\tilde u_{[2]} \\
\vdots \\
\tilde u_{[\ell-1]} \\
\tilde u_{[\ell]}
\end{bmatrix}
=
\begin{bmatrix}
f_{[1]} \\
f_{[2]} \\
\vdots \\
f_{[\ell-1]} \\
f_{[\ell]}
\end{bmatrix},
\end{gather*}
where $H_{[i,j]}$'s are the corresponding sparse blocks. Note that we use the bracket subscripts $[ \cdot ]$ to emphasize that the corresponding unknowns are grouped together in each slice.

\begin{figure}
[!ht]
\centering

\begin{tikzpicture}
[scale = 0.25]

\foreach \x in {-4,...,23}
\foreach \y in {-4,...,23}
\fill(\x,\y)[fill=blue!50]circle(0.25);

\foreach \x in {0,...,19}
\foreach \y in {0,...,19}
\fill(\x,\y)[fill=green!50]circle(0.25);

\foreach \x in {3.5,7.5,...,15.5}
\draw[dashed, line width = 1pt] (\x,-5.5) -- (\x,24.5);

\foreach \x in {1,...,5}
\draw(\x * 4 -2.5, -4.5)node[below]{$\cD_{\x}$};

\end{tikzpicture}
\caption{The grid points are divided into $5$ slices along the $x_1$ direction. Each of the middle slices contains $4$ layers, and each of the two boundary ones contains $4$ more attenuating layers in the PML region.}
\label{fig:2D_slice}
\end{figure}

We introduce the Schur complement $S_{[i]}$ and its inverse $T_{[i]}$ slice by slice recursively
\begin{gather*}
S_{[1]} = H_{[1,1]}, \quad T_{[1]} = S_{[1]}^{-1},\\
S_{[i]} = H_{[i,i]} - H_{[i,i-1]} T_{[i-1]} H_{[i-1,i]}, \quad T_{[i]} = S_{[i]}^{-1},\quad \tfor i=2,\dots, \ell.
\end{gather*}
Then we can solve $\tilde u$ by the Gaussian elimination
\begin{gather*}
\tilde u_{[1]} = T_{[1]} f_{[1]},\\
\tilde u_{[i]} = T_{[i]} (f_{[i]} - H_{[i,i-1]} \tilde u_{[i-1]}),\quad \tfor i = 2,\dots, \ell,\\
\tilde u_{[i]} = \tilde u_{[i]} - T_{[i]} (H_{[i,i+1]} \tilde u_{[i+1]}),\quad \tfor i= \ell-1,\dots, 1.
\end{gather*}

The expensive part of the above process is to compute $T_{[i]}$ and apply it to the vectors on
$\cD_i$. If say we formed $T_{[i]}$ directly, the computation would take $O(b^3n^3)$ steps and the
application $O(b^2n^2)$ steps. The sweeping factorization reduces the cost by approximating
$T_{[i]}$ with a subproblem. To introduce the approximation, we first make a key observation of the
operator $T_{[i]}$: inverting the top left $i\times i$ block of $H$, one notices that $T_{[i]}$
appears at the bottom right block of the resulting matrix. In other words
\begin{gather*}
H_{[1:i,1:i]}^{-1}=
\begin{bmatrix}
H_{[1,1]} & H_{[1,2]} \\
H_{[2,1]} & H_{[2,2]} & \ddots \\
 & \ddots & \ddots & \ddots \\
 & & \ddots & \ddots & H_{[i-1,i]} \\
 & & & H_{[i,i-1]} & H_{[i,i]}
\end{bmatrix}^{-1}
=
\begin{bmatrix}
* & * & \ldots & * & * \\
* & * & \ldots & * & * \\
\vdots & \vdots & \ddots & \vdots & \vdots \\
* & * & \ldots & * & * \\
* & * & \ldots & * & T_{[i]} 
\end{bmatrix}.
\end{gather*}
This means $T_{[i]}$ is the restriction of $H_{[1:i,1:i]}^{-1}$ to $\cD_i$. Think of $T_{[i]}$ as an operator from some input vector $v$ to $T_{[i]} v$ on the grid $\cD_i$. Then given $v$, we can compute $T_{[i]} v$ by solving the equation
\begin{equation}
\label{eqn:Hwv}
\begin{bmatrix}
H_{[1,1]} & H_{[1,2]} \\
H_{[2,1]} & H_{[2,2]} & \ddots \\
 & \ddots & \ddots & \ddots \\
 & & \ddots & \ddots & H_{[i-1,i]} \\
 & & & H_{[i,i-1]} & H_{[i,i]}
\end{bmatrix}
\begin{bmatrix}
* \\
* \\
\vdots \\
* \\
w
\end{bmatrix}
=
\begin{bmatrix}
0 \\
0 \\
\vdots \\
0 \\
v
\end{bmatrix}
\end{equation}
where $w$ is exactly equal to $T_{[i]} v$. That is to say, given $v$, we can find $T_{[i]} v$ by
padding $v$ with zeros on $\cD_{1:(i-1)}$, solving the unknowns on $\cD_{1:i}$ by \eqref{eqn:Hwv} and
then extracting the solution on $\cD_i$.

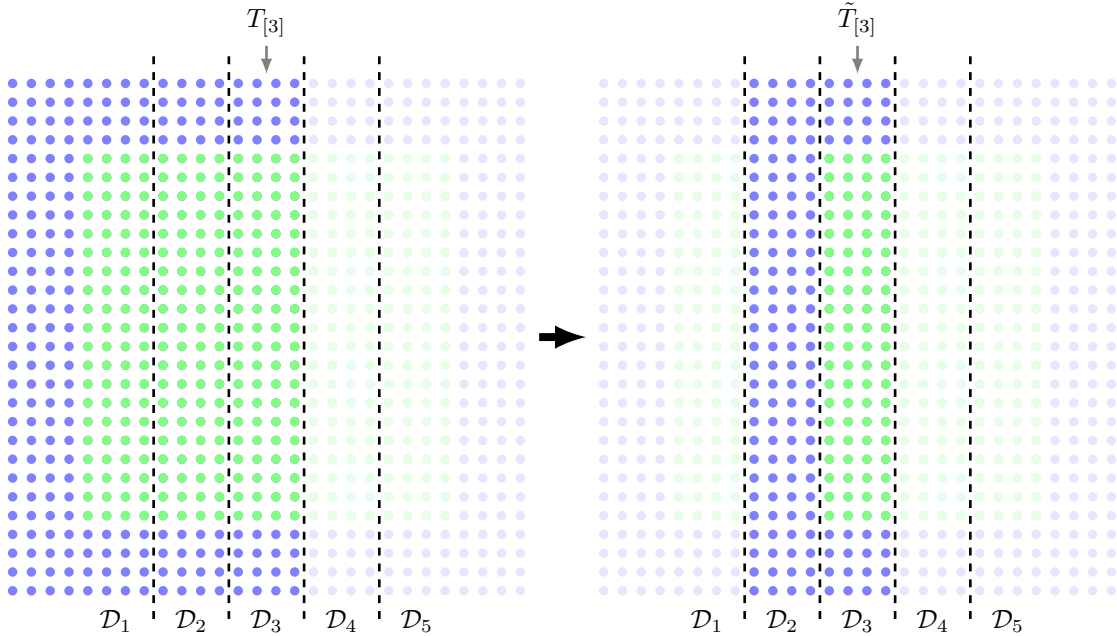
\begin{figure}
[!ht]
\centering

\begin{tikzpicture}
[scale = 0.25]

\foreach \x in {-4,...,23}
\foreach \y in {-4,...,23}
\fill(\x,\y)[fill=blue!10]circle(0.25);

\foreach \x in {0,...,19}
\foreach \y in {0,...,19}
\fill(\x,\y)[fill=green!10]circle(0.25);

\foreach \x in {-4,...,11}
\foreach \y in {-4,...,23}
\fill(\x,\y)[fill=blue!50]circle(0.25);

\foreach \x in {0,...,11}
\foreach \y in {0,...,19}
\fill(\x,\y)[fill=green!50]circle(0.25);

\foreach \x in {3.5,7.5,...,15.5}
\draw[dashed, line width = 1pt] (\x,-5.5) -- (\x,24.5);

\foreach \x in {1,...,5}
\draw(\x * 4 -2.5, -4.5)node[below]{$\cD_{\x}$};
\draw(9.5,25)node[above]{$T_{[3]}$}[->,>=latex,gray,line width=1pt] -- (9.5,23.5);

\draw[->, >=latex, line width = 3pt](24,9.5)--(26.5,9.5);
\end{tikzpicture}
\begin{tikzpicture}
[scale = 0.25]

\foreach \x in {-4,...,23}
\foreach \y in {-4,...,23}
\fill(\x,\y)[fill=blue!10]circle(0.25);

\foreach \x in {0,...,19}
\foreach \y in {0,...,19}
\fill(\x,\y)[fill=green!10]circle(0.25);

\foreach \x in {4,...,11}
\foreach \y in {-4,...,23}
\fill(\x,\y)[fill=blue!50]circle(0.25);

\foreach \x in {8,...,11}
\foreach \y in {0,...,19}
\fill(\x,\y)[fill=green!50]circle(0.25);

\foreach \x in {3.5,7.5,...,15.5}
\draw[dashed, line width = 1pt] (\x,-5.5) -- (\x,24.5);

\foreach \x in {1,...,5}
\draw(\x * 4 -2.5, -4.5)node[below]{$\cD_{\x}$};

\draw(9.5,25)node[above]{${\tilde T}_{[3]}$}[->,>=latex,gray,line width=1pt] -- (9.5,23.5);

\end{tikzpicture}
\caption{This figure is an illustration of the moving PML method. Left: $T_{[3]}$ is the restriction
  of $H_{[1:3,1:3]}^{-1}$ to the discrete domain $\cD_3$. Right: We move the PML adjacent to $\cD_3$
  to induce the approximation operator $\tilde T_{[3]}$ as the restriction to $\cD_3$ of the
  corresponding subproblem on $\cD_{2:3}$.}
\label{fig:move_PML}
\end{figure}

Note that the right-hand side of \eqref{eqn:Hwv} is zero on $\cD_{1:(i-1)}$, thus the only role of
the first $i-1$ blocks of equations in \eqref{eqn:Hwv} is to induce the radiation condition at the
left boundary of $\cD_i$ implicitly. To simulate this radiation condition, one can directly put the
PML region to the left side of $\cD_i$ instead of putting it far away on $\cD_1$. That's the key
idea of the sweeping factorization: move the PML region adjacent to the domain of interest $\cD_i$
and approximate the operator $T_{[i]}$ by solving a much smaller system compared to \eqref{eqn:Hwv}
(see Figure \ref{fig:move_PML}).

By introducing the modified plain waves, we can build the local stencils for points in the auxiliary
PML region on the left of $\cD_i$ similar to what was done in Section \ref{subsubsec:PML_stencil}. A
subtle difference is that, the local spacial frequency is perturbed to $\omega \sqrt{1 - m(x)}$
instead of $\omega$ at location $x$, and we need to use this local frequency to build the local
stencil for each point.

To save computational cost of the stencil construction, we do not use the exact value of the local
frequency. Though building local stencil in the PML region with the exact local frequency takes only
constant steps per point in theory, the constant is not small since it involves finding the kernel
of a $8\times 9$ matrix. Instead, we consider the square frequency range:
\begin{gather*}
[\quad \omega^2(1 - \max\{m(x)\})\quad ,\quad  \omega^2(1 - \min\{m(x)\})\quad].
\end{gather*}
We choose some samples uniformly from this range interval, and build local stencils only for those
samples. Then for each point in the PML region, we assign the stencil to be the one from the samples
with the closest local square frequency value, a technique introduced earlier in
\cite{liu2016localized}. In practice, only $n$ samples will be enough for an accurate
approximation. So it only takes $O(b n)$ steps to build these stencils, which is negligible compared
to the problem size $O(n^2)$. An intuition of why we only need $n$ samples is that, $O(\omega^2 /
n)$ is the size of the variation in one neighborhood $\mu_i$ on average, so there's no need to make
the sampling scale smaller than that.

With the auxiliary PML region on the left of $\cD_i$, we can solve a much smaller system instead of
solving \eqref{eqn:Hwv}. In our setting, the set of the auxiliary PML points for $\cD_i$ is just
$\cD_{i-1}$ since the width of the PML region is the same as $\cD_{i-1}$. The auxiliary system can be
written as
\begin{equation}
\label{eqn:Hwv_approx}
\begin{bmatrix}
\tilde H_{[i-1,i-1]} & \tilde H_{[i-1,i]}\\
H_{[i,i-1]} & H_{[i,i]}
\end{bmatrix}
\begin{bmatrix}
*\\
w
\end{bmatrix}
=
\begin{bmatrix}
0\\
v
\end{bmatrix},
\end{equation}
where the bottom block of equations is inherited from \eqref{eqn:Hwv}, and the top block is defined by the local stencils of the second type in the auxiliary PML region, the role of which is to simulate the radiation boundary condition on the left of $\cD_i$.

A minor problem here is that the auxiliary PML region for $\cD_2$ consists only the normal layers in $\cD_1$ rather than all the layers, so \eqref{eqn:Hwv_approx} needs a slight modification for $i = 2$: we restrict the columns of $H_{[2,1]}$ to the normal layers in $\cD_1$ so that the two blocks are compatible. This problem is inessential and the patch here is only to make the discussion strictly correct. In practice, the width of the slices and the PML regions can be rather flexible.

Equation \eqref{eqn:Hwv_approx} defines an approximating operator $\tilde T_{[i]}: v \to w$ for
$i\in 1\dots, \ell$ by restricting the system \eqref{eqn:Hwv_approx} on $\cD_{(i-1):i}$ to
$\cD_i$. Note that for $i = 1$, $\tilde T_{[1]}$ is exactly equal to $T_{[1]}$ if we treat $\cD_0$
as $\emptyset$ naturally. Compared to \eqref{eqn:Hwv}, Equation \eqref{eqn:Hwv_approx} is a much
smaller quasi-1D problem, which can be solved efficiently with the LU factorization.

\subsection{Putting together}
We now have all the tools needed to design a linear-complexity preconditioner for the discretized
Lippmann-Schwinger equation \eqref{eqn:LS_d}. The setup and application processes of the
preconditioner are given by Algorithms \ref{alg:2D_setup} and \ref{alg:2D_app} respectively. The
slice width $b$ is typically a small integer less than $10$, thus both the setup and the
application costs are linear.

\begin{algorithm}
[!ht]
  \caption{Setup of the preconditioner for the system \eqref{eqn:LS_d}. Complexity $=O(b^2N)$. }
  \label{alg:2D_setup}
  \begin{algorithmic}[1]
  \STATE
  Compute the stencils $\alpha$ and $\beta$. Complexity $=O(N)$.
  \STATE
  Compute the PML stencils $\gamma$ for different local frequency samples and different positions of complex stretching. Complexity $=O(b n)$.
  \STATE
  Divide the domain into $\ell$ slices as $\cD_1\dots,\cD_\ell$.
  \STATE
  Define the approximating operator $\tilde T_{[i]}$ by the sweeping factorization below from Step \ref{algstep:setup_start} to \ref{algstep:setup_end}. Complexity $=O(b^2 N)$:
  
  \FOR {$i = 1,\dots,\ell$}
  \label{algstep:setup_start}
    \STATE
    Pad $\cD_i$ with auxiliary PML points in $\cD_{i-1}$ to form the subproblem \eqref{eqn:Hwv_approx} where the auxiliary PML stencils are built with samples closet to the local square frequency values.
    \STATE
    Compute the LU factorization of \eqref{eqn:Hwv_approx}, which defines the solution operator. Restricting the solution operator to $\cD_i$ induces $\tilde T_{[i]}$.
    \ENDFOR
  \label{algstep:setup_end}
  \end{algorithmic}
\end{algorithm}

\begin{algorithm}
[!ht]
  \caption{Application of the preconditioner. Complexity $=O(b N)$.}
  \label{alg:2D_app}
  \begin{algorithmic}[1]
    \STATE
    Form the right-hand side $f$ of \eqref{eqn:Huf}. Complexity $=O(N)$.
    \STATE
    Solve the linear system \eqref{eqn:Huf} approximately by the process below from Step \ref{algstep:app_start} to \ref{algstep:app_end}.\\
    Complexity $=O(b N)$.
    \STATE
  \label{algstep:app_start}
    $\tilde u_{[1]} = \tilde T_{[1]} f_{[1]}$
    \FOR {$i=2,\dots,\ell$}
    \STATE
    $\tilde u_{[i]} = \tilde T_{[i]} (f_{[i]} - H_{[i, i-1]} \tilde u_{[i-1]})$
    \ENDFOR
    \FOR{$i=\ell-1,\dots,1$}
    \STATE
    $\tilde u_{[i]} = \tilde u_{[i]} - \tilde T_{[i]}(H_{[i,i+1]} \tilde u_{[i+1]})$
    \ENDFOR
  \label{algstep:app_end}
    \STATE Now $\tilde u$ is an approximation to $H^{-1} f$. Extract the solution of $\tilde u$ on $\cD$ as the output.
  \end{algorithmic}
\end{algorithm}

We would like to make some comments below for the actual implementation of the algorithm.
\begin{enumerate}
\item
The algorithm presented above constructs the sweeping factorization along the $x_1$ direction from
left to right. Indeed, since we have radiation conditions on all sides of the domain, we can
construct the factorization from both sides and sweep toward the middle slice. The two sweeping
fronts can be processed independently until they meet in the middle, where they exchange some local
information in the middle slice and then sweep back to the boundaries independently. This is
potentially helpful for the parallelization of the algorithm.
\item
The widths of the slices and the auxiliary PML regions are completely arbitrary. There are two
reasons why we set them to be $b$ uniformly in what above. The first is for the simplicity of
discussion. The second is that, given the PML width $b$, it is optimal to set the width of each
slice to be also $b$ to minimize the setup and application costs of the preconditioner. In practice,
it may not be possible to uniformly divide the domain where each slice contains $b$ layers
exactly. In that case, we change the widths of one or two slices accordingly, which has negligible
effect to the cost and efficiency of the preconditioner.
\item
The constructions of the stencils $\alpha, \beta$ and $\gamma$, though depending on $n$ and
$\omega$, are essentially independent of the velocity field $c(x)$. First, the computation of
$\alpha$ and $\beta$ only involves the free space Green's function $G(x)$, where the velocity field
is completely irrelevant. Next, for the local PML stencils $\gamma$, they might depend on $c(x)$
slightly, but only on the range as we see from the sampling process. In practice $c(x)=\Theta(1)$,
so the range is actually bounded for fixed $\omega$. Thus we can precompute the stencils without
given the velocity field. This means that the stencil construction only needs a fixed cost for given
problem size, which can be eliminated from the setup process of the algorithm for the input $c(x)$.
\end{enumerate}

\subsection{Numerical results}
\label{subsec:2D_numerical}

In this section we present the numerical results in 2D. The algorithm is implemented in MATLAB and the tests are performed on a 2.4-GHz server. We force MATLAB to use only one computational thread to test the sequential time cost. The preconditioner is combined with the standard GMRES solver with relative tolerance $10^{-6}$ and restart value $20$. The domain is discretized with $h = \lambda / 8$ where $\lambda = 2\pi / \omega$ is the typical wavelength.

We choose $b = 8$ as the width of the slices and the auxiliary PML regions. This corresponds to about one wavelength width for the PML
regions and the slices used in the sweeping preconditioner. The sweeping factorization is built with
two fronts sweeping toward the middle slice, and the middle slice is padded with auxiliary PMLs on
both sides for the corresponding quasi-1D subproblem.

Four velocity fields are tested in 2D, which are
\begin{enumerate}
[(i)]
\item
  \label{LS:vel:2D_1}
  A converging Gaussian centered at $(0.5,0.5)$.
\item
  \label{LS:vel:2D_2}
  A diverging Gaussian centered at $(0.5,0.5)$.
\item
  \label{LS:vel:2D_3}
  $32$ randomly placed converging Gaussians with narrow width.
\item
  \label{LS:vel:2D_4}
  A random velocity field that is equal to $1$ at $\partial\Omega$.
\end{enumerate}
The incoming wave $u_I(x)$ for each test is a plane wave shooting downward at frequency
$\omega$. The test results are given in Tables \ref{LS:tab:2D_1}, \ref{LS:tab:2D_2}, \ref{LS:tab:2D_3} and \ref{LS:tab:2D_4} respectively. The notations in the tables are listed as follows.
\begin{itemize}
\item
$\omega$ is the angular frequency.
\item
$N$ is the number of unknowns.
\item
$T_{\text{setup}}$ is the setup cost of the preconditioner in seconds.
\item
$T_{\text{apply}}$ is the application cost of the preconditioner in seconds.
\item
$N_{\text{iter}}$ is the iteration number.
\item
$T_{\text{solve}}$ is the solve cost of the preconditioner in seconds.
\end{itemize}

\begin{table}
[!ht]

\centering
\includegraphics[width=0.38\textwidth]{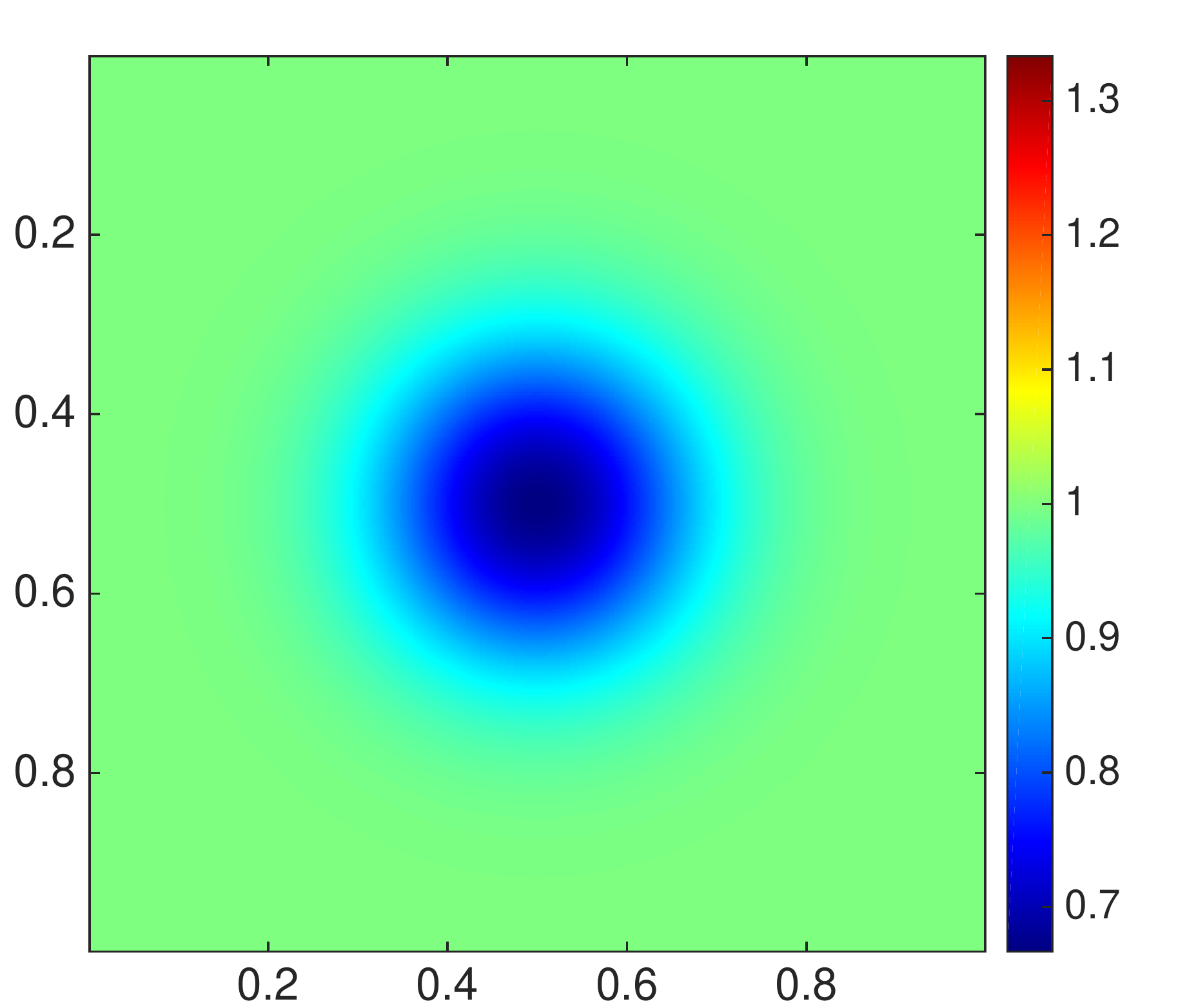}
\includegraphics[width=0.38\textwidth]{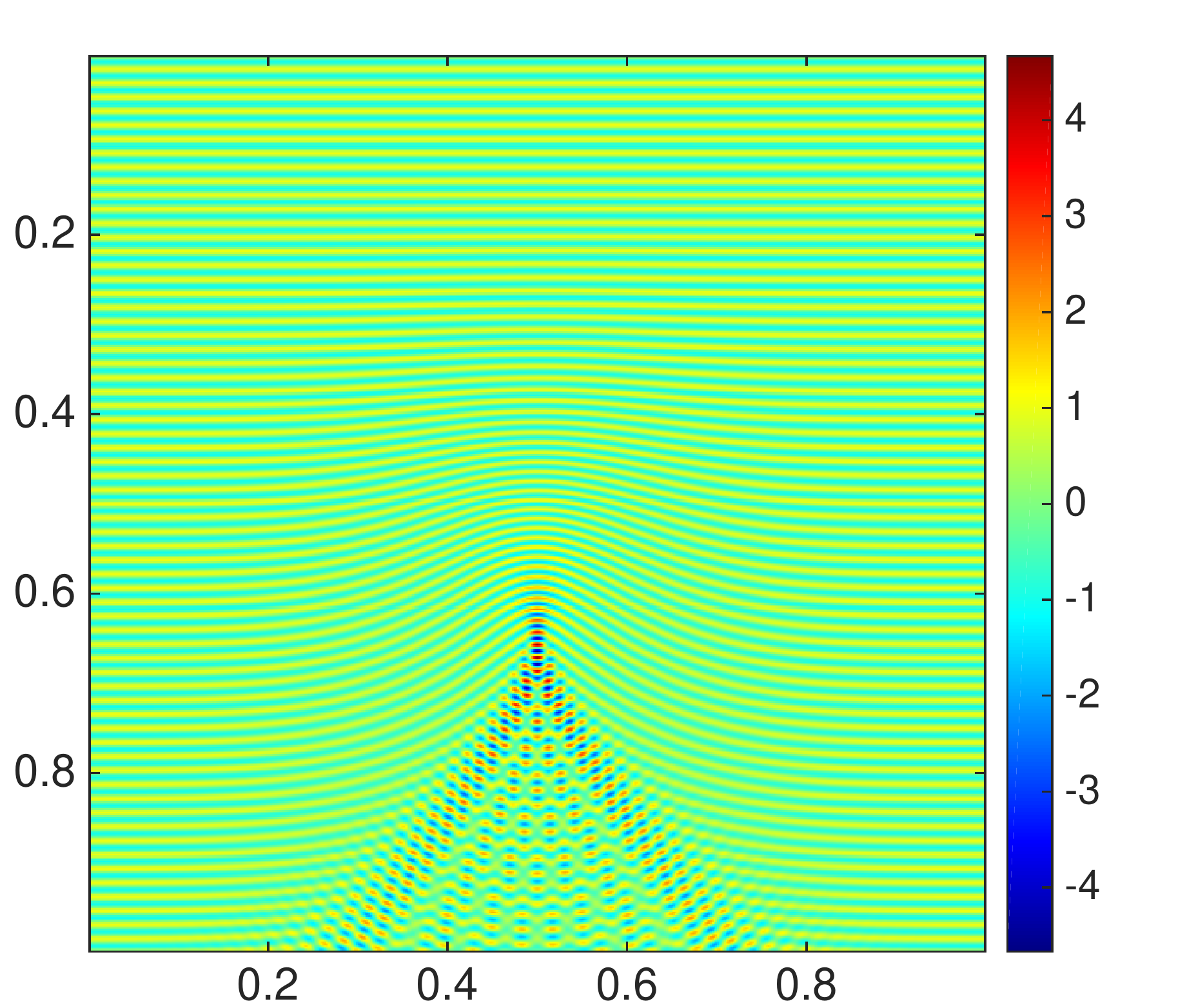}

\begin{tabular}{cc | cc | cc}
\hline
\hline
$\omega/(2\pi)$ & $N$ & $T_{\text{setup}}$ & $T_{\text{apply}}$ & $N_{\text{iter}}$ & $T_{\text{solve}}$ \\
\hline
$16$ & $127^2$ & $1.03$e$-01$ & $1.47$e$-02$ & $5$ & $1.96$e$-01$ \\
$32$ & $255^2$ & $3.01$e$-01$ & $4.85$e$-02$ & $5$ & $4.27$e$-01$ \\
$64$ & $511^2$ & $1.23$e$+00$ & $1.82$e$-01$ & $5$ & $1.89$e$+00$ \\
$128$ & $1023^2$ & $4.80$e$+00$ & $7.29$e$-01$ & $6$ & $8.58$e$+00$ \\
$256$ & $2047^2$ & $1.95$e$+01$ & $2.90$e$+00$ & $7$ & $4.69$e$+01$ \\
\hline
\hline
\end{tabular}

\caption[Numerical results for velocity field \eqref{LS:vel:2D_1} in 2D.]{Numerical results for velocity field \eqref{LS:vel:2D_1} in 2D. Top: The velocity field $c(x)$ (left) and the total wave field $u(x)+u_I(x)$ (right) for $\omega/(2\pi) = 64$. Bottom: Table of the numerical results for different problem sizes.}
\label{LS:tab:2D_1}
\end{table}

\begin{table}
[!ht]

\centering
\includegraphics[width=0.38\textwidth]{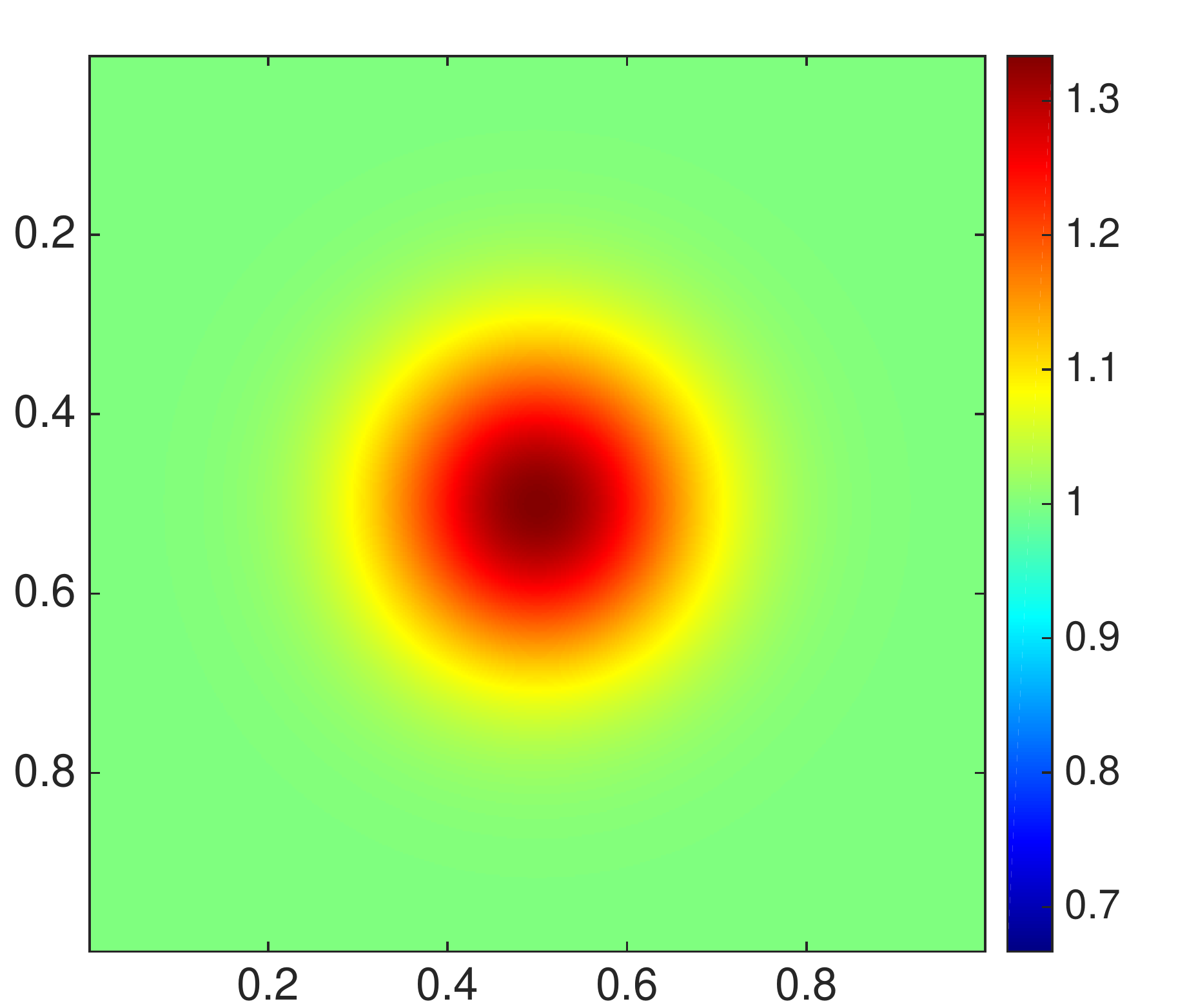}
\includegraphics[width=0.38\textwidth]{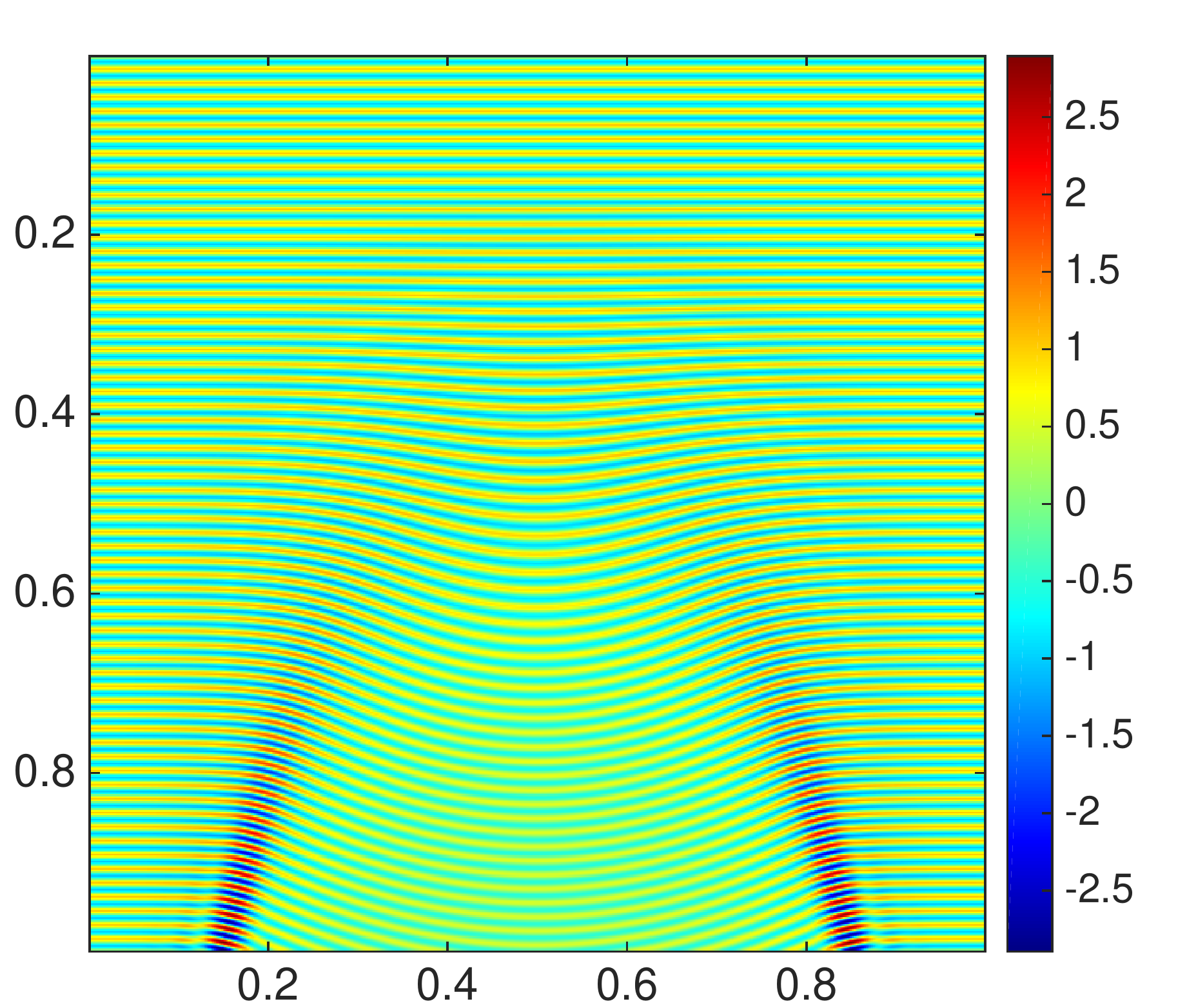}

\begin{tabular}{cc | cc | cc}
\hline
\hline
$\omega/(2\pi)$ & $N$ & $T_{\text{setup}}$ & $T_{\text{apply}}$ & $N_{\text{iter}}$ & $T_{\text{solve}}$ \\
\hline
$16$ & $127^2$ & $1.59$e$-01$ & $2.73$e$-02$ & $4$ & $1.79$e$-01$ \\
$32$ & $255^2$ & $5.39$e$-01$ & $4.86$e$-02$ & $4$ & $3.29$e$-01$ \\
$64$ & $511^2$ & $1.23$e$+00$ & $1.82$e$-01$ & $5$ & $1.61$e$+00$ \\
$128$ & $1023^2$ & $4.81$e$+00$ & $7.04$e$-01$ & $5$ & $7.19$e$+00$ \\
$256$ & $2047^2$ & $1.95$e$+01$ & $2.89$e$+00$ & $6$ & $4.01$e$+01$ \\
\hline
\hline
\end{tabular}

\caption[Numerical results for velocity field \eqref{LS:vel:2D_2} in 2D.]{Numerical results for velocity field \eqref{LS:vel:2D_2} in 2D. Top: The velocity field $c(x)$ (left) and the total wave field $u(x)+u_I(x)$ (right) for $\omega/(2\pi) = 64$. Bottom: Table of the numerical results for different problem sizes.}
\label{LS:tab:2D_2}
\end{table}

\begin{table}
[!ht]

\centering
\includegraphics[width=0.38\textwidth]{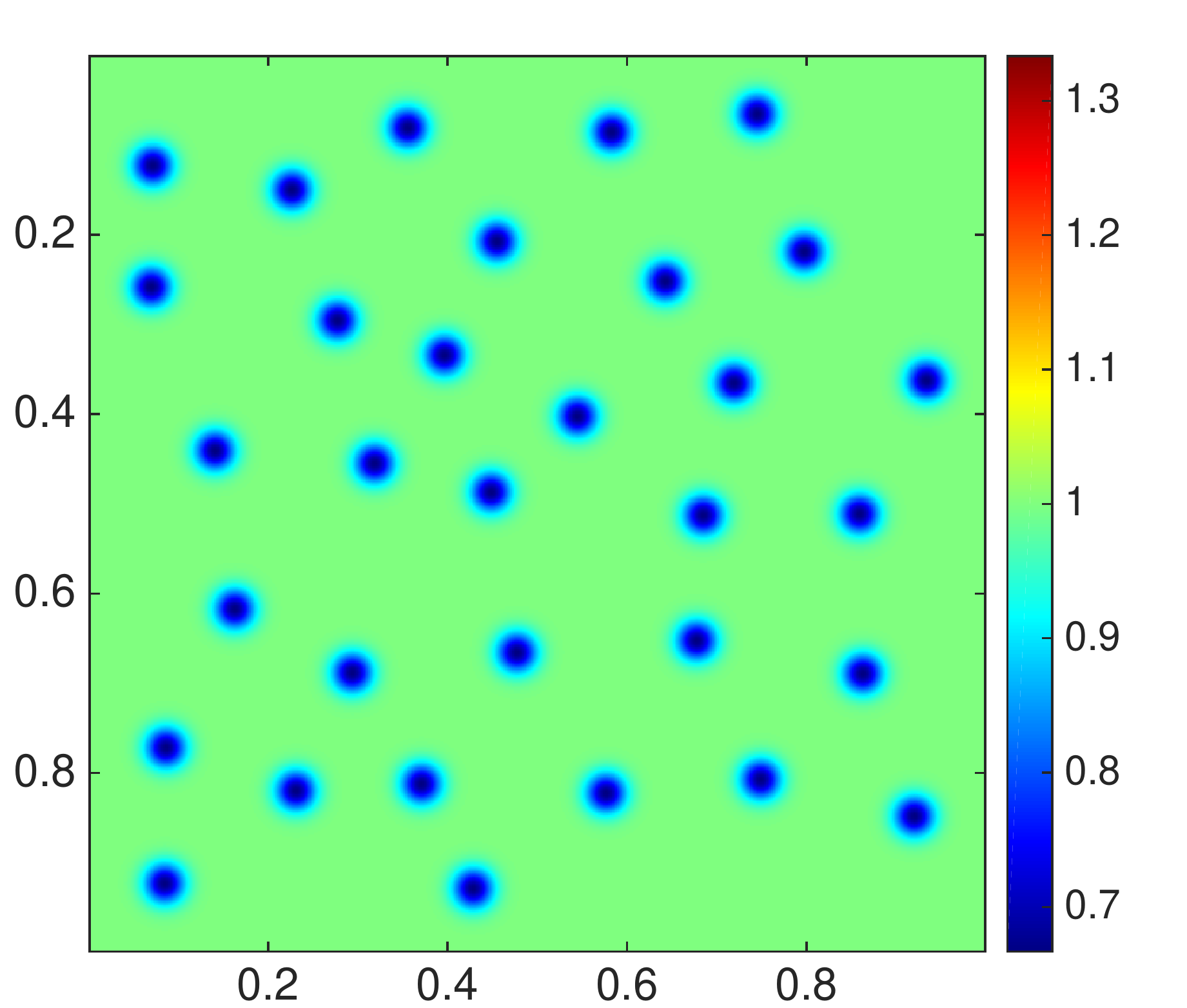}
\includegraphics[width=0.38\textwidth]{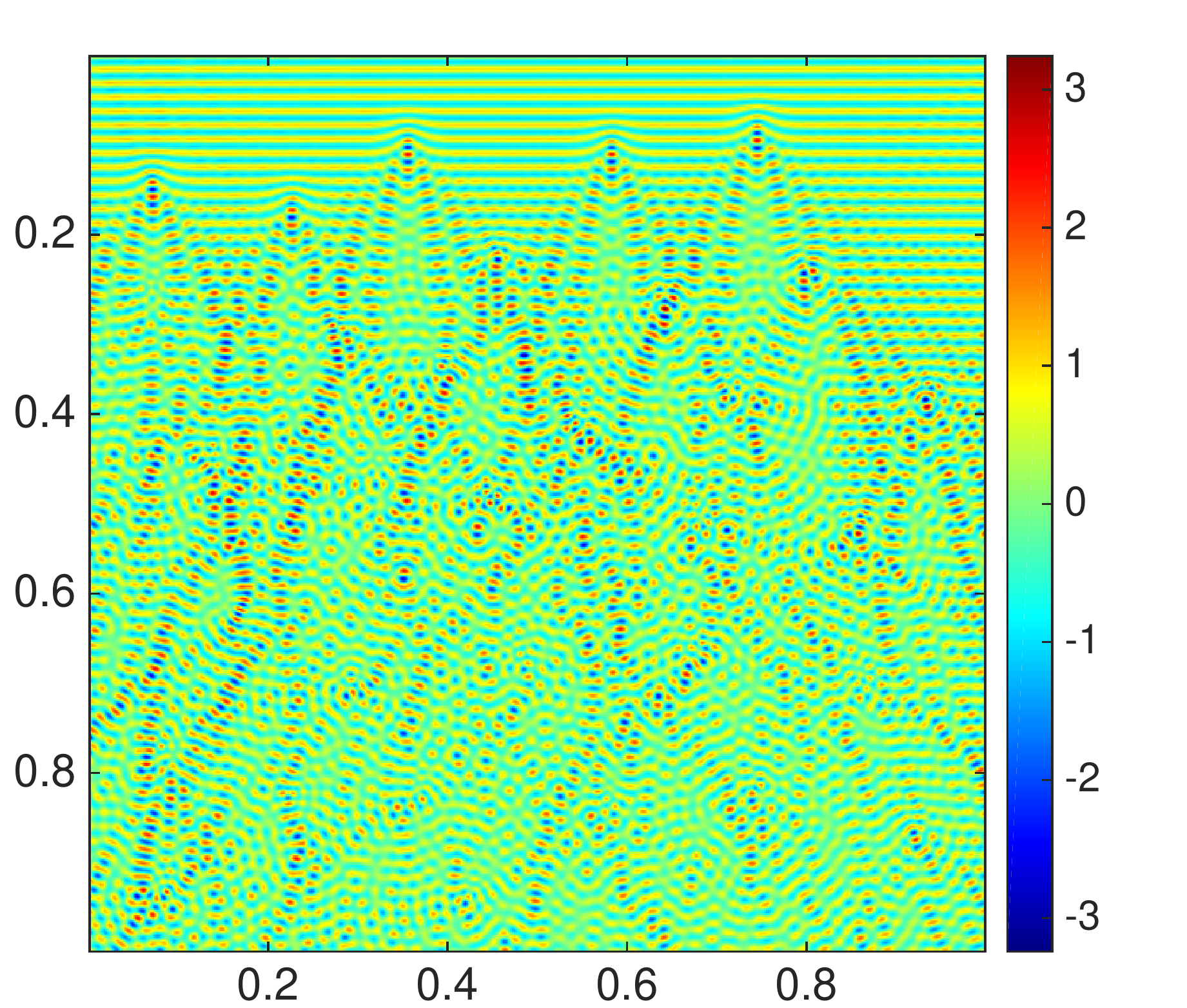}

\begin{tabular}{cc | cc | cc}
\hline
\hline
$\omega/(2\pi)$ & $N$ & $T_{\text{setup}}$ & $T_{\text{apply}}$ & $N_{\text{iter}}$ & $T_{\text{solve}}$ \\
\hline
$16$ & $127^2$ & $9.27$e$-02$ & $1.39$e$-02$ & $9$ & $2.02$e$-01$ \\
$32$ & $255^2$ & $3.01$e$-01$ & $4.88$e$-02$ & $8$ & $6.65$e$-01$ \\
$64$ & $511^2$ & $1.23$e$+00$ & $1.82$e$-01$ & $9$ & $2.95$e$+00$ \\
$128$ & $1023^2$ & $4.80$e$+00$ & $7.08$e$-01$ & $10$ & $1.43$e$+01$ \\
$256$ & $2047^2$ & $1.94$e$+01$ & $2.96$e$+00$ & $11$ & $7.57$e$+01$ \\
\hline
\hline
\end{tabular}

\caption[Numerical results for velocity field \eqref{LS:vel:2D_3} in 2D.]{Numerical results for velocity field \eqref{LS:vel:2D_3} in 2D. Top: The velocity field $c(x)$ (left) and the total wave field $u(x)+u_I(x)$ (right) for $\omega/(2\pi) = 64$. Bottom: Table of the numerical results for different problem sizes.}
\label{LS:tab:2D_3}
\end{table}

\begin{table}
[!ht]

\centering
\includegraphics[width=0.38\textwidth]{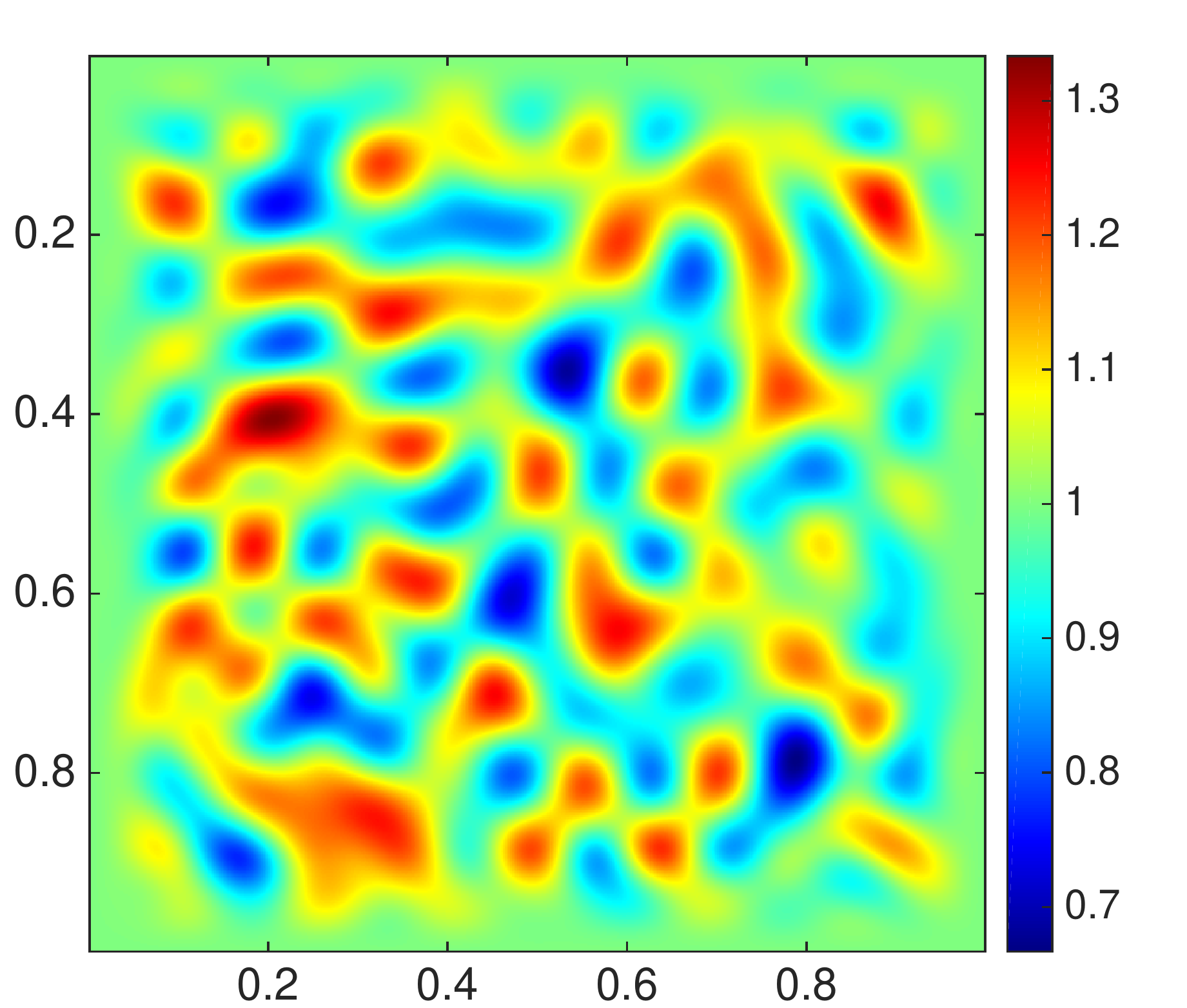}
\includegraphics[width=0.38\textwidth]{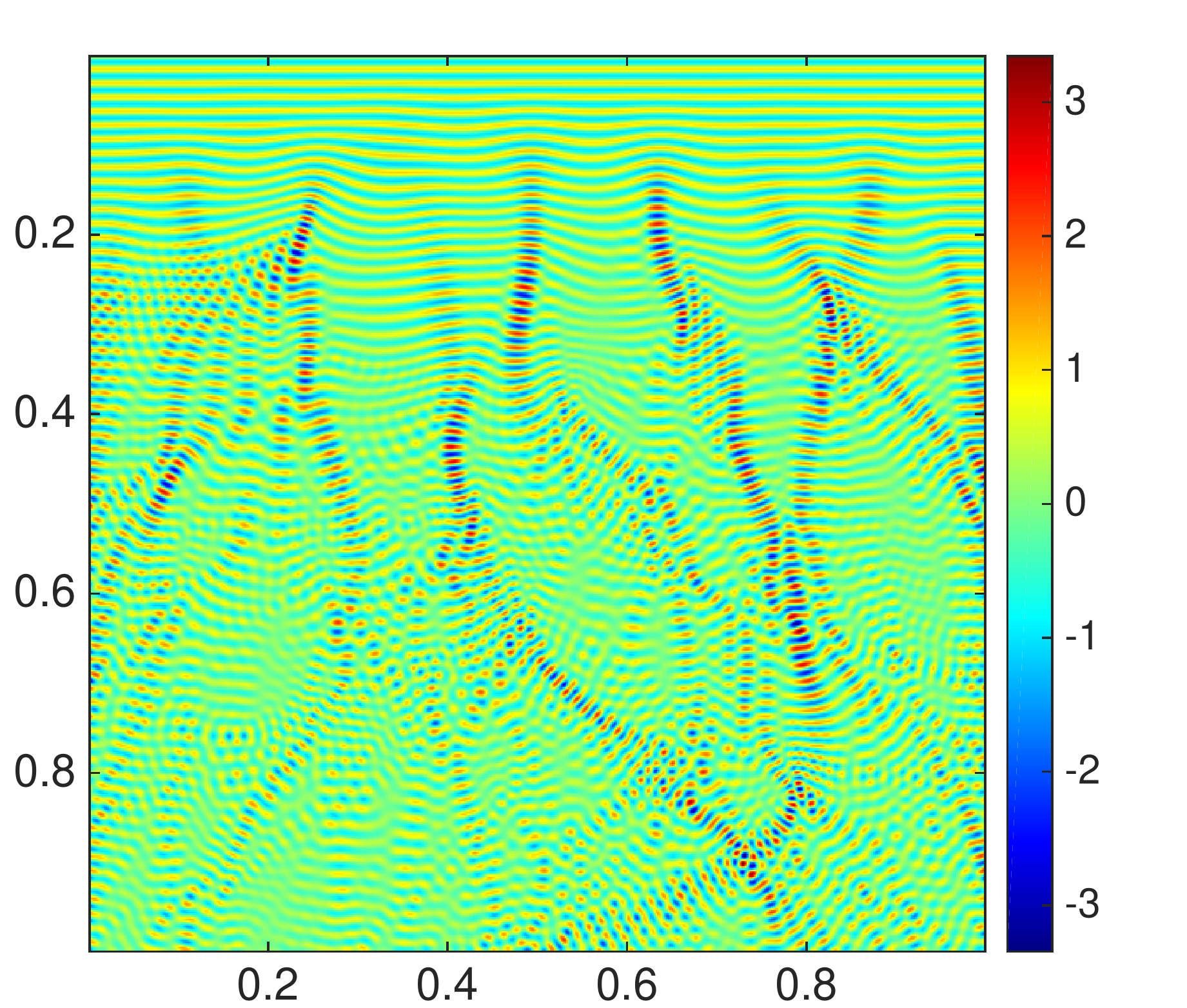}

\begin{tabular}{cc | cc | cc}
\hline
\hline
$\omega/(2\pi)$ & $N$ & $T_{\text{setup}}$ & $T_{\text{apply}}$ & $N_{\text{iter}}$ & $T_{\text{solve}}$ \\
\hline
$16$ & $127^2$ & $9.05$e$-02$ & $1.39$e$-02$ & $7$ & $1.54$e$-01$ \\
$32$ & $255^2$ & $2.98$e$-01$ & $4.88$e$-02$ & $7$ & $5.80$e$-01$ \\
$64$ & $511^2$ & $1.23$e$+00$ & $1.82$e$-01$ & $8$ & $2.65$e$+00$ \\
$128$ & $1023^2$ & $4.83$e$+00$ & $7.07$e$-01$ & $9$ & $1.26$e$+01$ \\
$256$ & $2047^2$ & $1.95$e$+01$ & $2.90$e$+00$ & $9$ & $6.09$e$+01$ \\
\hline
\hline
\end{tabular}

\caption[Numerical results for velocity field \eqref{LS:vel:2D_4} in 2D.]{Numerical results for velocity field \eqref{LS:vel:2D_4} in 2D. Top: The velocity field $c(x)$ (left) and the total wave field $u(x)+u_I(x)$ (right) for $\omega/(2\pi) = 64$. Bottom: Table of the numerical results for different problem sizes.}
\label{LS:tab:2D_4}
\end{table}

From the numerical tests we observe that both the setup time and the application time scale linearly
in $N$, which are in accordance with the complexity analyses. More importantly, the iteration numbers change only slightly as the problem size grows, almost independent of $\omega$.

We notice that the iteration number also depends on the velocity field. For simple fields such as
the diverging Gaussian, it requires less iterations compared to more complicated fields such as the
narrower converging Gaussians. This makes sense intuitively since converging lenses and velocity
fields with drastic local variations increase the oscillations and refractions of the wave field,
thus the corresponding systems are harder to solve. In addition, for the sweeping factorization to
work well, we need to assume that there are no strong reflections and refractions during the
transmission of the waves so that the auxiliary PMLs in the intermediate slices can make correct
approximations to the true underlying DtN maps. In practice, moderate amount of wave-ray bendings
can be taken care of by a few more iterations as we see in the tests for the multiple diverging
Gaussians and the random field. If the velocity field is even worse, for example, if the field has
large region of strong discontinuities, then neither will the Nystr\"om method be able to give an
accurate discretization scheme, nor can the sweeping factorization provide an accurate approximating
solution due to the strong reflections caused by the discontinuities. Thus for our preconditioner to
work, we require certain smoothness from the velocity fields. Nonetheless, as we can tell from the
numerical examples, the preconditioner works well even when the fields have drastic transitions in
narrow regions. So this approach can be widely applied to many use cases.

\section{Preconditioner in 3D}
\label{sec:3D}

This section presents the preconditioner in 3D. As we see from Section \ref{sec:2D}, the approach is
essentially dimension independent and it can be easily generalized to 3D. We will keep the
description short, mainly emphasizing the differences compared to the 2D case so that the reader can
get the central idea effortlessly. The 3D numerical results for both the recursive approach and the
non-recursive approach of the sweeping factorization are provided in the second part of this
section.

\subsection{Problem formulation, sparsification and sweeping factorization}
In this section we formulate the approach in 3D. All the notations in 2D can be easily reused
without causing any ambiguities. We will keep them unless otherwise stated.

We assume $\Omega=(0,1)^3$ contains the support of $m(x)$. The domain is discretized with step size
$h = O(1/\omega)$ in each dimension. A similar quadrature correction formula is used for the central
weight of the Green's function, which gives an accuracy of $O(h^4)$.

For the sparsification process, the first type of stencils $\alpha$ and $\beta$ can be constructed
similarly, where now each neighborhood $\mu_i$ has $27$ points. For the second type of stencils in
the PML region, we use the modified plain waves in 3D, defined similarly as
\begin{gather*}
F^\sigma(x) \coloneqq \exp(\ii \omega (r \cdot x^\sigma)), \quad \|r\|_2 = 1,
\end{gather*}
where now $x^\sigma$ and $r$ are in $\bR^3$, and $x^\sigma$ is stretched to the complex plane from
$x$ for all three coordinates. In 2D, the stencil $\gamma$ is defined as the kernel vector which
annihilates the independent waves shooting toward the eight most common directions. This can be done
similarly in 3D. We now need a set of $26$ directions, which is defined as
\begin{gather*}
R \coloneqq \left\{\dfrac{(r_1,r_2,r_3)}{\sqrt{r_1^2+r_2^2+r_3^2}} : (r_1,r_2,r_3) \in \{-1,0, 1\}^3 \setminus \{(0, 0, 0)\}\right\}.
\end{gather*}
In other words, these are the directions shooting from the center of a neighborhood to the $26$
boundary neighbor points.

The computational cost of constructing the stencils in 3D seems higher due to more degrees of
freedom and larger size of the neighborhoods. But indeed, the relative cost compared to the
sweeping factorization is lower than the 2D case, let alone that the stencil computations are
independent of the velocity field and they can be done by a once-in-a-life-time preprocessing.

For the sweeping factorization, the domain are now divided into $\ell$ quasi-2D slices. The auxiliary
PMLs are padded to each slice similarly. Each subproblem is quasi-2D, which can be solved
efficiently by the nested dissection algorithm with $O(b^3n^3)$ setup cost and $O(b^2n^2\log n)$
application cost. Consisting of $\ell \approx n / b$ subproblems, the whole process has a total setup
cost $O(b^2n^4) = O(b^2 N^{4/3})$ and application cost $O(b n^3 \log n) = O(b N \log N)$. Note that
the direct use of the nested dissection algorithm to the 3D sparse system costs $O(N^2)$ for setup
and $O(N^{4/3})$ for solve. The sweeping factorization drastically reduces the costs by dimension
reduction.

For each of the quasi-2D problem, we can sweep similarly along the $x_2$ direction, reducing it to
$\ell$ quasi-1D subproblems. This reduces the setup cost to $O(b^4 N)$ and the application cost to
$O(b^2 N)$, which are both linear in $N$, but more sensitive to the slice width $b$. We call this
the {\em recursive} approach \cite{Liu2015} while the one in the previous paragraph as {\em
  non-recursive}.

\subsection{Numerical results}
\label{subsec:3D_numerical}

In this section we present the numerical results in 3D. The test configurations are the same as Section \ref{subsec:2D_numerical} unless otherwise stated. In the 3D tests, we set $b = 4$ for the slice width and PML width.

The four velocity fields tested are
\begin{enumerate}
[(i)]
\item
\label{LS:vel:3D_1}
A converging Gaussian centered at $(0.5,0.5,0.5)$.
\item
\label{LS:vel:3D_2}
A diverging Gaussian centered at $(0.5,0.5,0.5)$.
\item
\label{LS:vel:3D_3}
$256$ randomly placed converging Gaussians of  narrow width.
\item
\label{LS:vel:3D_4}
A random velocity field that is equal to $1$ at $\partial\Omega$.
\end{enumerate}
The right-hand side is a plain wave shooting downward at frequency $\omega$.

The tests of the non-recursive approach are given in Tables \ref{LS:tab:3D_1},\ref{LS:tab:3D_2}, \ref{LS:tab:3D_3} and \ref{LS:tab:3D_4}, and the ones of the recursive approach are in Tables \ref{LS:tab:3D_recursive_1}, \ref{LS:tab:3D_recursive_2}, \ref{LS:tab:3D_recursive_3}, \ref{LS:tab:3D_recursive_4}, where the relative costs compared to the non-recursive approach are also
listed as percentages, together with the iteration numbers of the non-recursive ones in the
parentheses for the convenience of comparison.

\begin{table}
[!ht]

\centering
\includegraphics[width=0.38\textwidth]{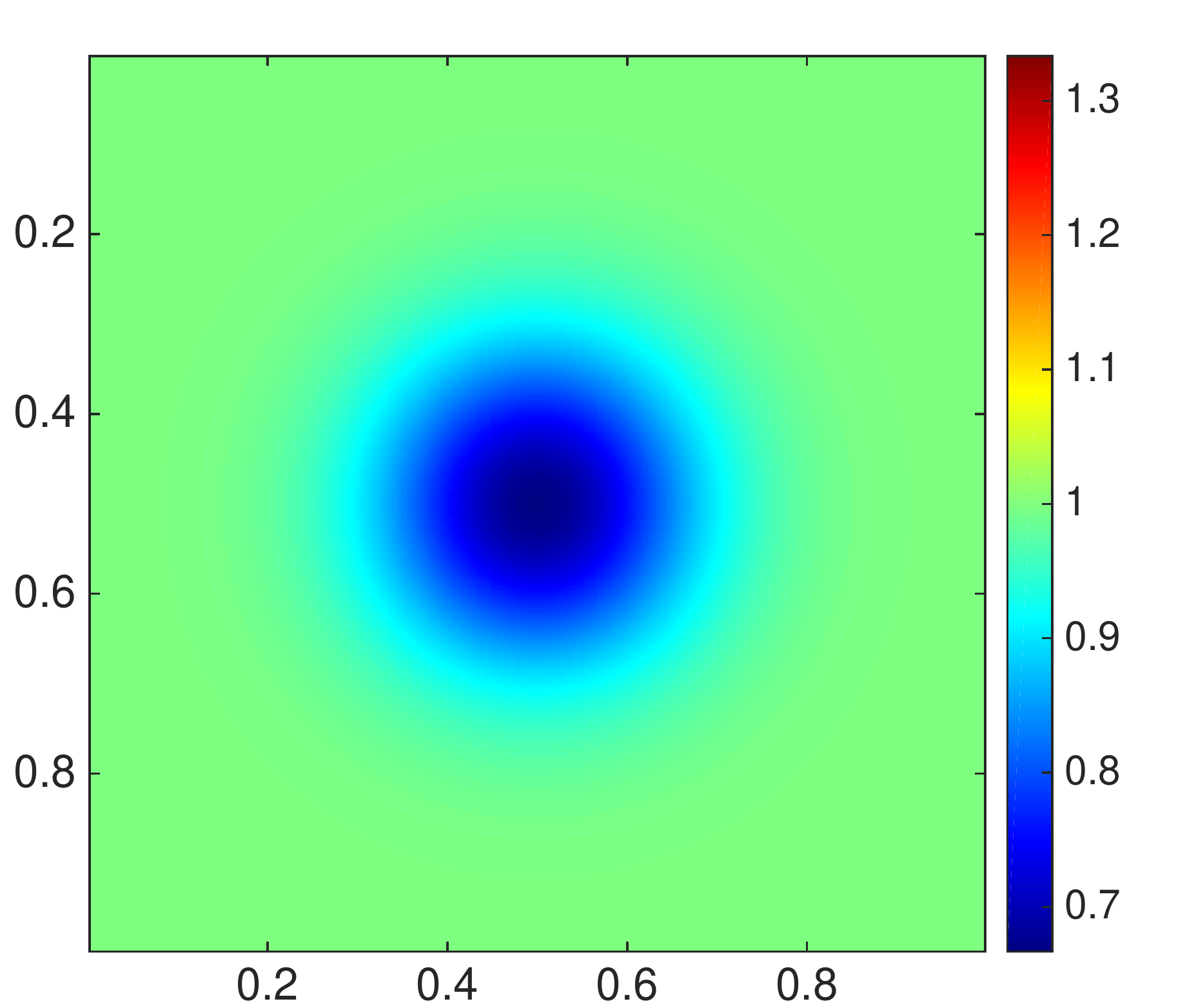}
\includegraphics[width=0.38\textwidth]{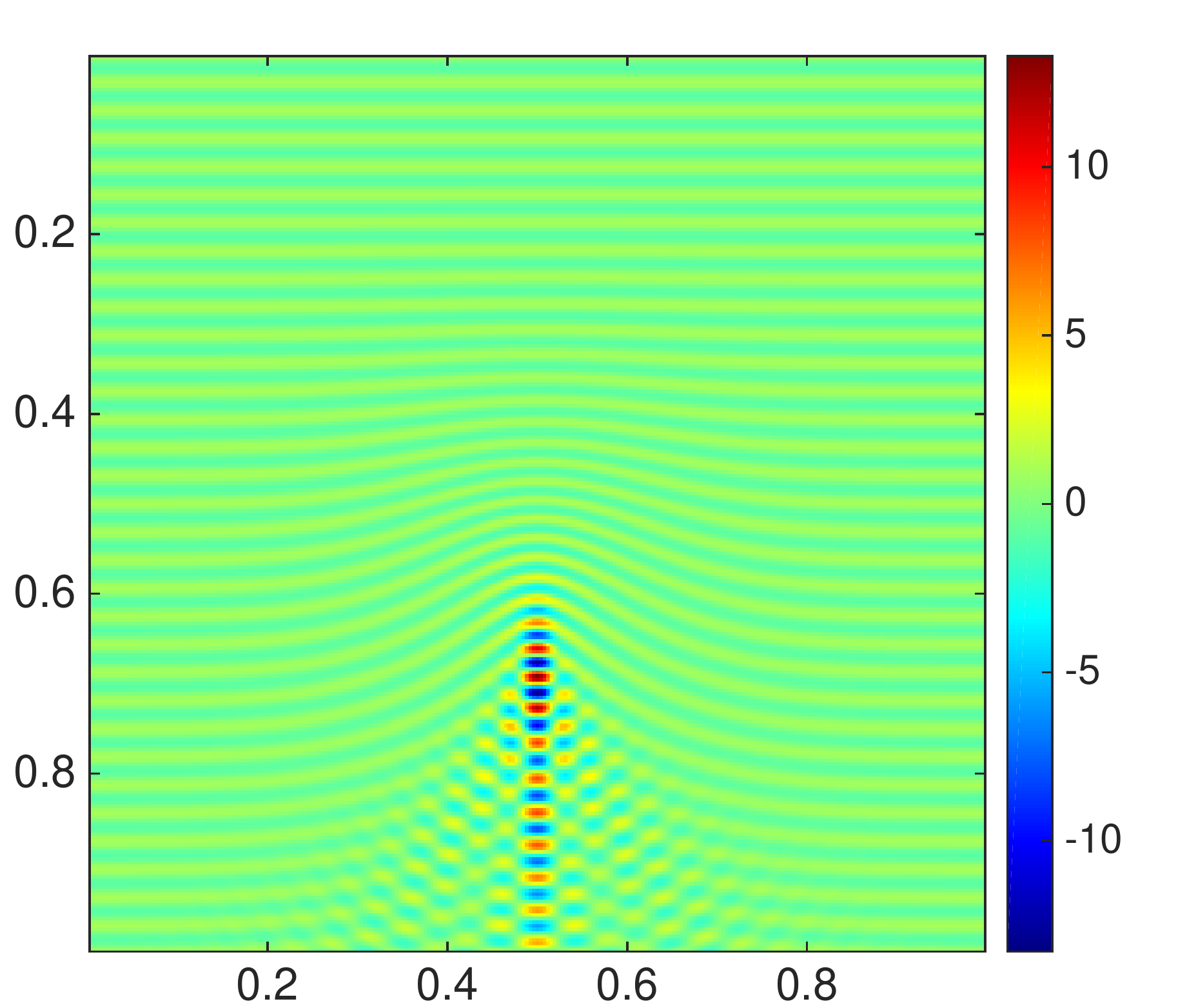}

\begin{tabular}{cc | cc | cc}
\hline
\hline
$\omega/(2\pi)$ & $N$ & $T_{\text{setup}}$ & $T_{\text{apply}}$ & $N_{\text{iter}}$ & $T_{\text{solve}}$ \\
\hline
$4$ & $31^3$ & $6.85$e$+00$ & $4.20$e$-01$ & $5$ & $2.46$e$+00$ \\
$8$ & $63^3$ & $5.74$e$+01$ & $3.28$e$+00$ & $5$ & $1.78$e$+01$ \\
$16$ & $127^3$ & $5.97$e$+02$ & $2.74$e$+01$ & $5$ & $1.57$e$+02$ \\
$32$ & $255^3$ & $7.24$e$+03$ & $2.49$e$+02$ & $6$ & $1.68$e$+03$ \\
\hline
\hline
\end{tabular}

\caption[Numerical results for velocity field \eqref{LS:vel:3D_1} in 3D with the non-recursive approach.]{Numerical results for velocity field \eqref{LS:vel:3D_1} in 3D with the non-recursive approach. Top: The velocity field $c(x)$ (left) and the total wave field $u(x)+u_I(x)$ (right) in a cross-section view for $\omega/(2\pi) = 32$. Bottom: Table of the numerical results for different problem sizes.}
\label{LS:tab:3D_1}
\end{table}

\begin{table}
[!ht]

\centering
\includegraphics[width=0.38\textwidth]{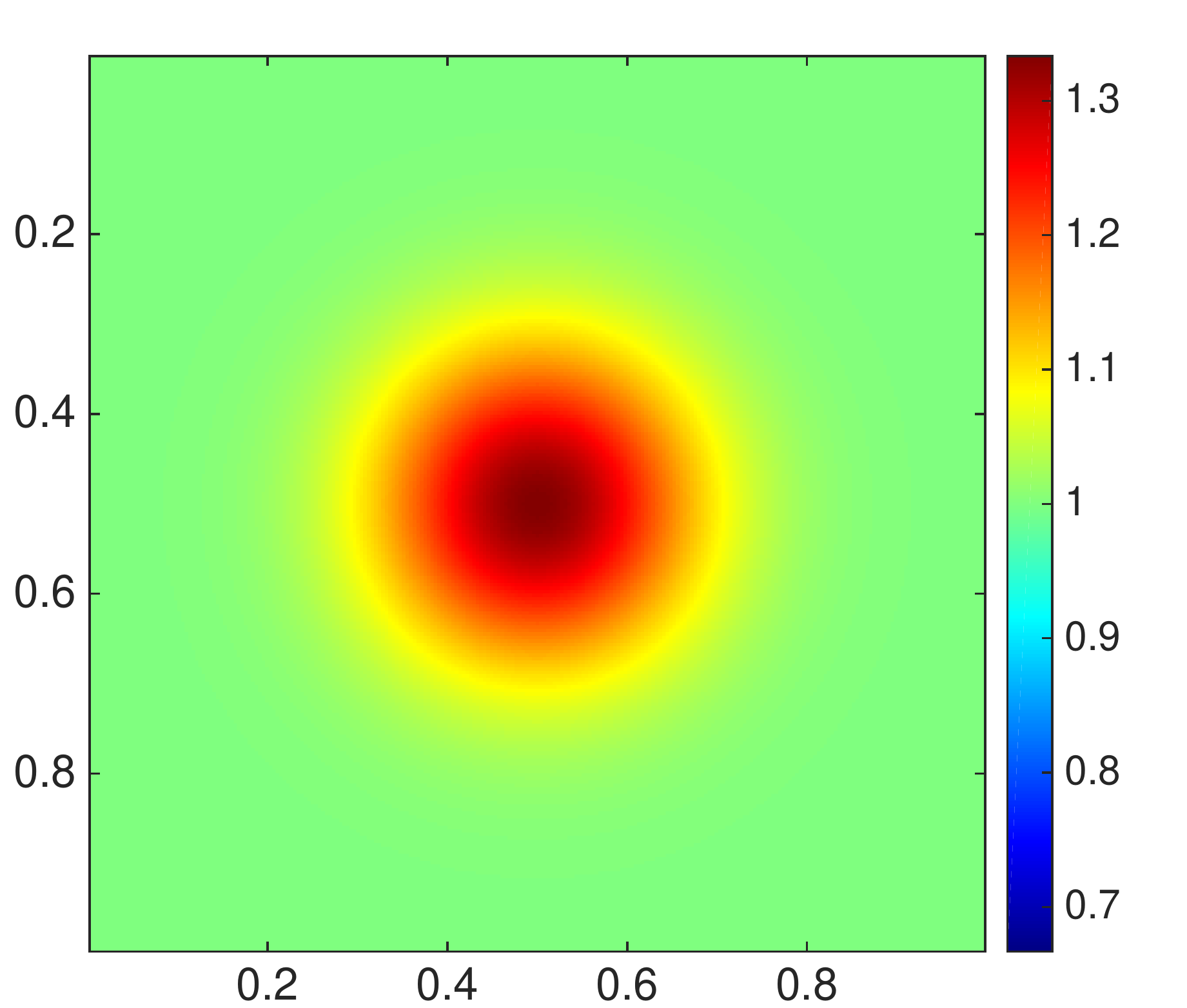}
\includegraphics[width=0.38\textwidth]{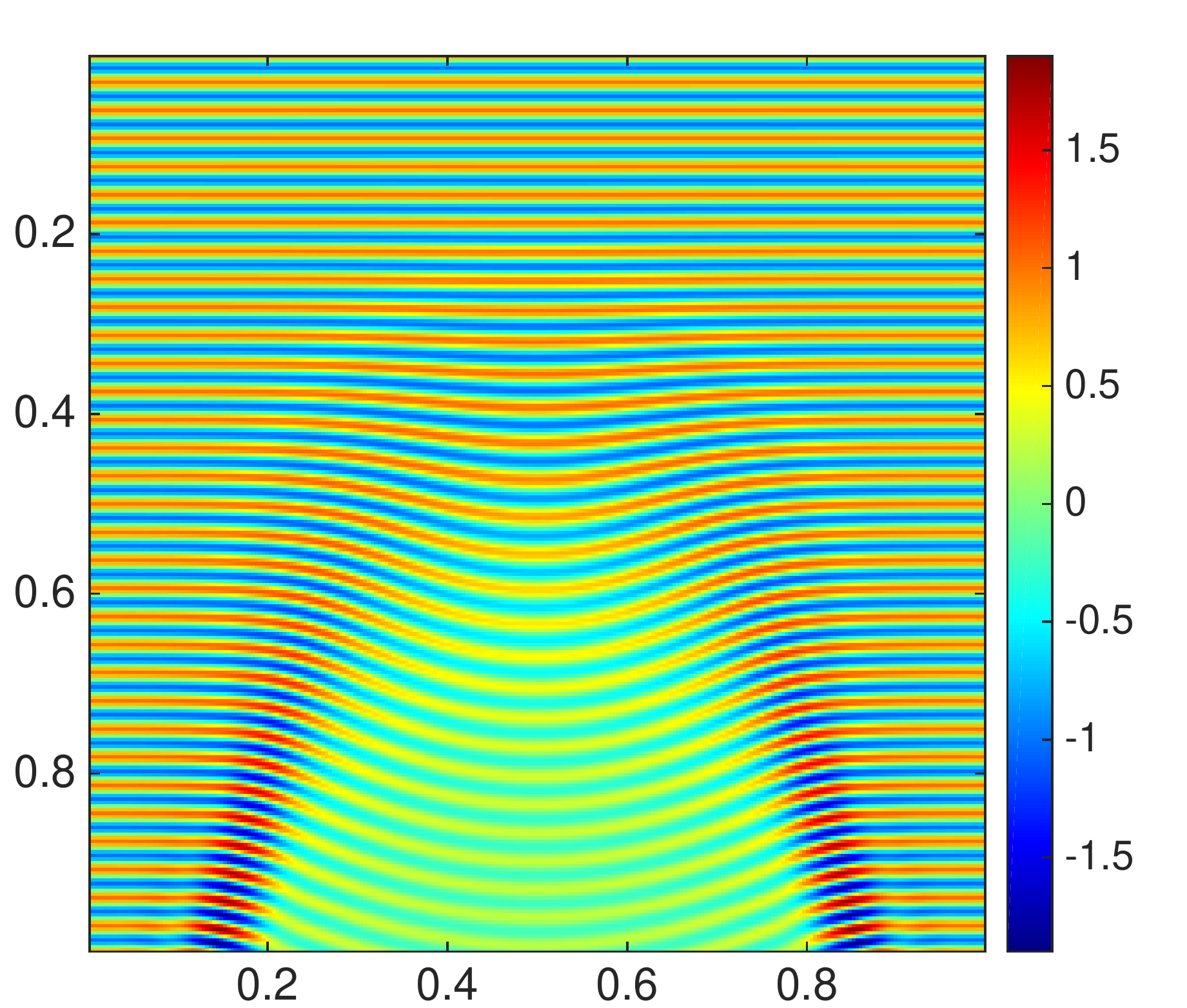}

\begin{tabular}{cc | cc | cc}
\hline
\hline
$\omega/(2\pi)$ & $N$ & $T_{\text{setup}}$ & $T_{\text{apply}}$ & $N_{\text{iter}}$ & $T_{\text{solve}}$ \\
\hline
$4$ & $31^3$ & $7.04$e$+00$ & $4.27$e$-01$ & $4$ & $1.85$e$+00$ \\
$8$ & $63^3$ & $5.85$e$+01$ & $3.26$e$+00$ & $5$ & $1.77$e$+01$ \\
$16$ & $127^3$ & $5.99$e$+02$ & $2.70$e$+01$ & $5$ & $1.54$e$+02$ \\
$32$ & $255^3$ & $7.24$e$+03$ & $2.49$e$+02$ & $5$ & $1.41$e$+03$ \\
\hline
\hline
\end{tabular}

\caption[Numerical results for velocity field \eqref{LS:vel:3D_2} in 3D with the non-recursive approach.]{Numerical results for velocity field \eqref{LS:vel:3D_2} in 3D with the non-recursive approach. Top: The velocity field $c(x)$ (left) and the total wave field $u(x)+u_I(x)$ (right) in a cross-section view for $\omega/(2\pi) = 32$. Bottom: Table of the numerical results for different problem sizes.}
\label{LS:tab:3D_2}
\end{table}

\begin{table}
[!ht]

\centering
\includegraphics[width=0.38\textwidth]{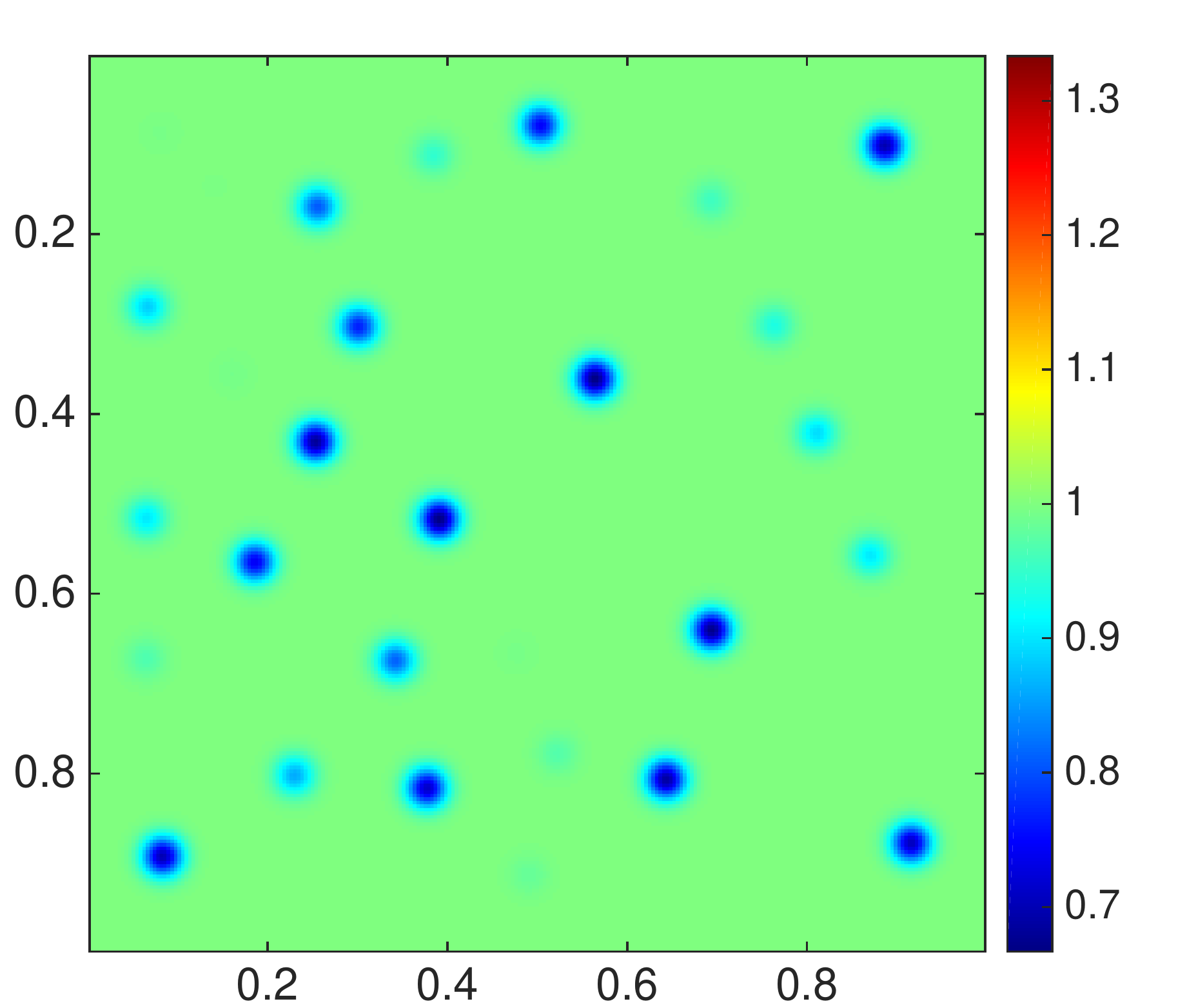}
\includegraphics[width=0.38\textwidth]{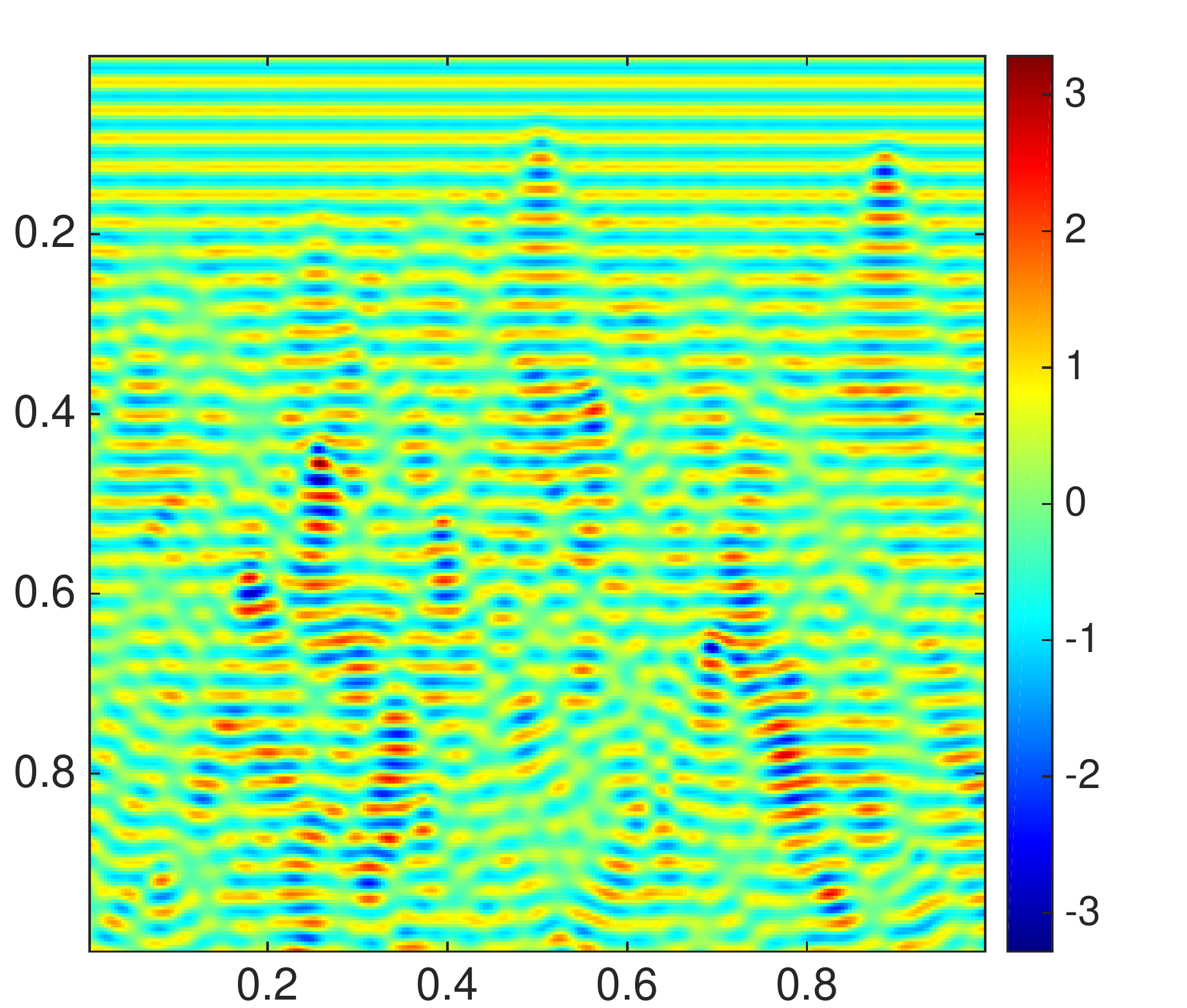}

\begin{tabular}{cc | cc | cc}
\hline
\hline
$\omega/(2\pi)$ & $N$ & $T_{\text{setup}}$ & $T_{\text{apply}}$ & $N_{\text{iter}}$ & $T_{\text{solve}}$ \\
\hline
$4$ & $31^3$ & $6.84$e$+00$ & $4.23$e$-01$ & $6$ & $2.92$e$+00$ \\
$8$ & $63^3$ & $5.71$e$+01$ & $3.31$e$+00$ & $8$ & $2.87$e$+01$ \\
$16$ & $127^3$ & $5.93$e$+02$ & $2.75$e$+01$ & $8$ & $2.50$e$+02$ \\
$32$ & $255^3$ & $7.19$e$+03$ & $2.48$e$+02$ & $8$ & $2.35$e$+03$ \\
\hline
\hline
\end{tabular}

\caption[Numerical results for velocity field \eqref{LS:vel:3D_3} in 3D with the non-recursive approach.]{Numerical results for velocity field \eqref{LS:vel:3D_3} in 3D with the non-recursive approach. Top: The velocity field $c(x)$ (left) and the total wave field $u(x)+u_I(x)$ (right) in a cross-section view for $\omega/(2\pi) = 32$. Bottom: Table of the numerical results for different problem sizes.}
\label{LS:tab:3D_3}
\end{table}

\begin{table}
[!ht]

\centering
\includegraphics[width=0.38\textwidth]{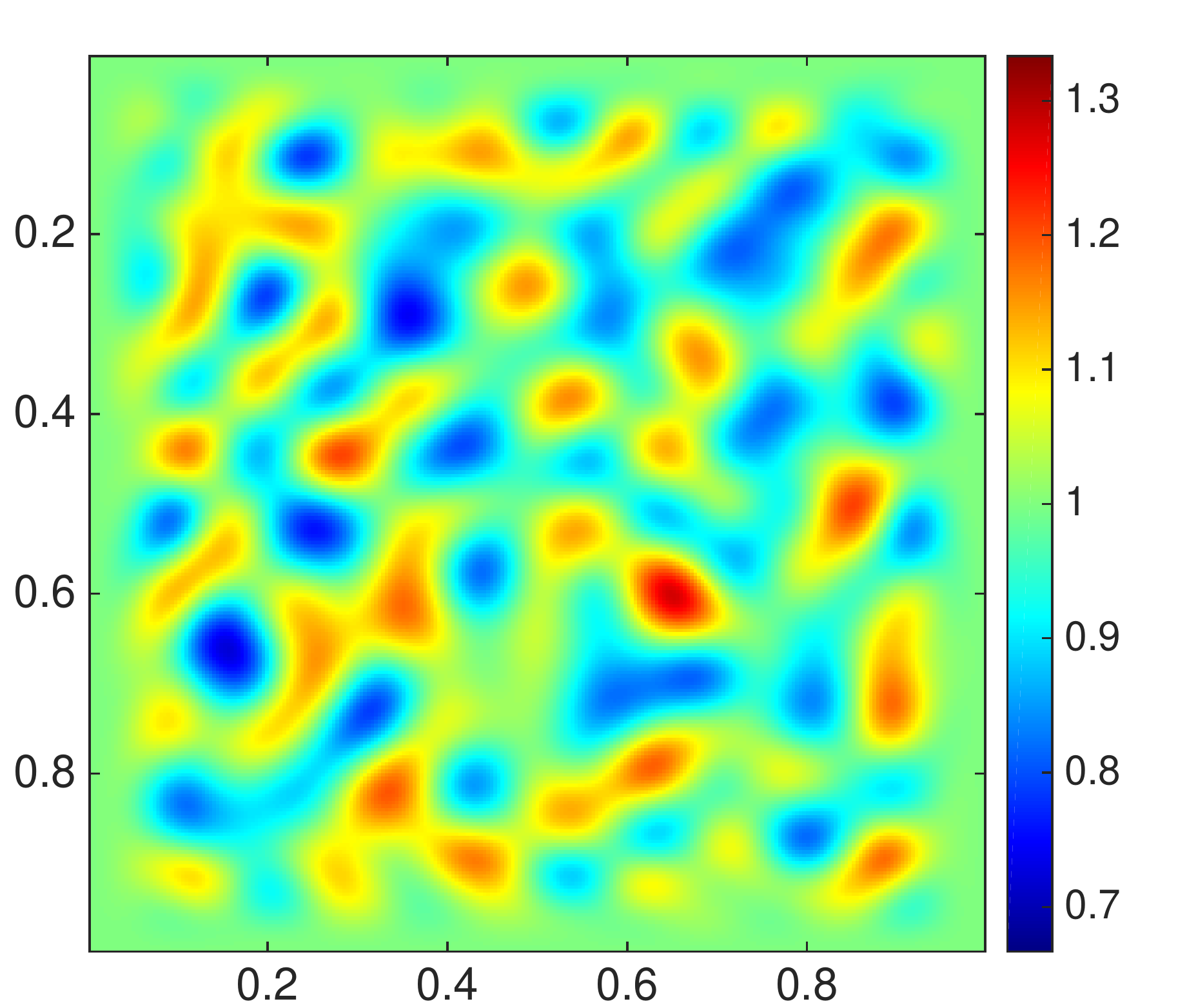}
\includegraphics[width=0.38\textwidth]{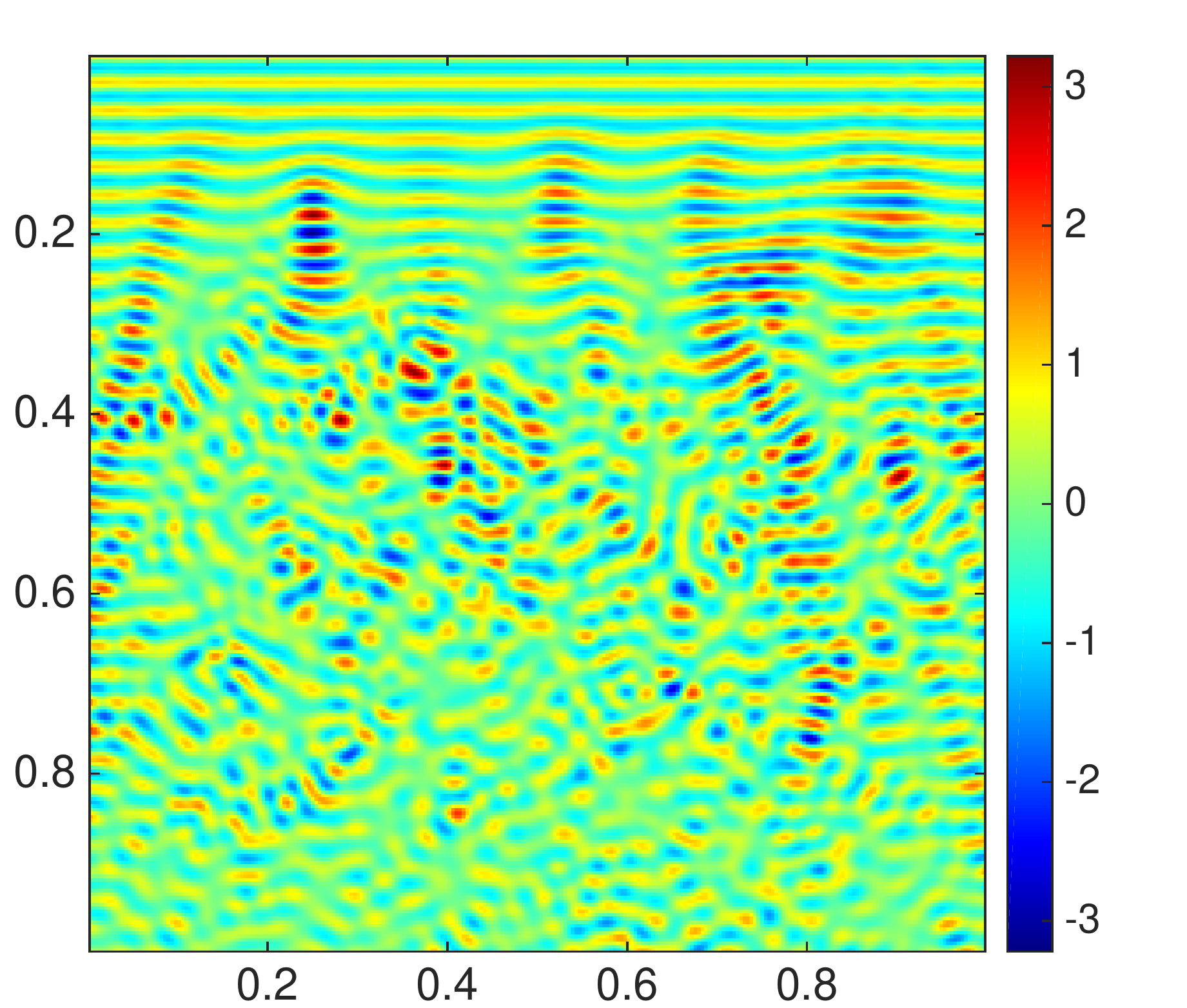}

\begin{tabular}{cc | cc | cc}
\hline
\hline
$\omega/(2\pi)$ & $N$ & $T_{\text{setup}}$ & $T_{\text{apply}}$ & $N_{\text{iter}}$ & $T_{\text{solve}}$ \\
\hline
$4$ & $31^3$ & $7.09$e$+00$ & $4.33$e$-01$ & $10$ & $4.64$e$+00$ \\
$8$ & $63^3$ & $5.85$e$+01$ & $3.32$e$+00$ & $10$ & $3.60$e$+01$ \\
$16$ & $127^3$ & $5.99$e$+02$ & $2.76$e$+01$ & $9$ & $2.81$e$+02$ \\
$32$ & $255^3$ & $7.20$e$+03$ & $2.48$e$+02$ & $9$ & $2.65$e$+03$ \\
\hline
\hline
\end{tabular}

\caption[Numerical results for velocity field \eqref{LS:vel:3D_4} in 3D with the non-recursive approach.]{Numerical results for velocity field \eqref{LS:vel:3D_4} in 3D with the non-recursive approach. Top: The velocity field $c(x)$ (left) and the total wave field $u(x)+u_I(x)$ (right) in a cross-section view for $\omega/(2\pi) = 32$. Bottom: Table of the numerical results for different problem sizes.}
\label{LS:tab:3D_4}
\end{table}

\begin{table}
[!ht]

\centering

\begin{tabular}{cc | cc | cc}
\hline
\hline
$\omega/(2\pi)$ & $N$ & $T_{\text{setup}}$ & $T_{\text{apply}}$ & $N_{\text{iter}}$ & $T_{\text{solve}}$ \\
\hline
$4$ & $31^3$ & $2.79$e$+00$ ($41$\%) & $3.44$e$-01$ ($82$\%) & $5$ ($5$) & $2.01$e$+00$ ($82$\%) \\
$8$ & $63^3$ & $1.62$e$+01$ ($28$\%) & $2.25$e$+00$ ($69$\%) & $5$ ($5$) & $1.29$e$+01$ ($73$\%) \\
$16$ & $127^3$ & $1.10$e$+02$ ($18$\%) & $1.64$e$+01$ ($60$\%) & $5$ ($5$) & $1.01$e$+02$ ($65$\%) \\
$32$ & $255^3$ & $8.23$e$+02$ ($11$\%) & $1.24$e$+02$ ($50$\%) & $6$ ($6$) & $9.33$e$+02$ ($56$\%) \\
\hline
\hline
\end{tabular}

\caption[Numerical results for velocity field \eqref{LS:vel:3D_1} in 3D with the recursive approach.]{Numerical results for velocity field \eqref{LS:vel:3D_1} in 3D with the recursive approach.}
\label{LS:tab:3D_recursive_1}
\end{table}

\begin{table}
[!ht]

\centering

\begin{tabular}{cc | cc | cc}
\hline
\hline
$\omega/(2\pi)$ & $N$ & $T_{\text{setup}}$ & $T_{\text{apply}}$ & $N_{\text{iter}}$ & $T_{\text{solve}}$ \\
\hline
$4$ & $31^3$ & $2.81$e$+00$ ($40$\%) & $3.38$e$-01$ ($79$\%) & $5$ ($4$) & $1.84$e$+00$ ($100$\%) \\
$8$ & $63^3$ & $1.63$e$+01$ ($28$\%) & $2.27$e$+00$ ($70$\%) & $5$ ($5$) & $1.28$e$+01$ ($72$\%) \\
$16$ & $127^3$ & $1.11$e$+02$ ($19$\%) & $1.64$e$+01$ ($61$\%) & $5$ ($5$) & $1.02$e$+02$ ($66$\%) \\
$32$ & $255^3$ & $8.18$e$+02$ ($11$\%) & $1.24$e$+02$ ($50$\%) & $6$ ($5$) & $9.31$e$+02$ ($66$\%) \\
\hline
\hline
\end{tabular}

\caption[Numerical results for velocity field \eqref{LS:vel:3D_2} in 3D with the recursive approach.]{Numerical results for velocity field \eqref{LS:vel:3D_2} in 3D with the recursive approach.}
\label{LS:tab:3D_recursive_2}
\end{table}

\begin{table}
[!ht]

\centering

\begin{tabular}{cc | cc | cc}
\hline
\hline
$\omega/(2\pi)$ & $N$ & $T_{\text{setup}}$ & $T_{\text{apply}}$ & $N_{\text{iter}}$ & $T_{\text{solve}}$ \\
\hline
$4$ & $31^3$ & $2.80$e$+00$ ($41$\%) & $3.36$e$-01$ ($79$\%) & $6$ ($6$) & $2.19$e$+00$ ($75$\%) \\
$8$ & $63^3$ & $1.63$e$+01$ ($29$\%) & $2.25$e$+00$ ($68$\%) & $8$ ($8$) & $2.02$e$+01$ ($70$\%) \\
$16$ & $127^3$ & $1.11$e$+02$ ($19$\%) & $1.64$e$+01$ ($60$\%) & $8$ ($8$) & $1.62$e$+02$ ($65$\%) \\
$32$ & $255^3$ & $8.20$e$+02$ ($11$\%) & $1.24$e$+02$ ($50$\%) & $8$ ($8$) & $1.24$e$+03$ ($53$\%) \\
\hline
\hline
\end{tabular}

\caption[Numerical results for velocity field \eqref{LS:vel:3D_3} in 3D with the recursive approach.]{Numerical results for velocity field \eqref{LS:vel:3D_3} in 3D with the recursive approach.}
\label{LS:tab:3D_recursive_3}
\end{table}

\begin{table}
[!ht]

\centering

\begin{tabular}{cc | cc | cc}
\hline
\hline
$\omega/(2\pi)$ & $N$ & $T_{\text{setup}}$ & $T_{\text{apply}}$ & $N_{\text{iter}}$ & $T_{\text{solve}}$ \\
\hline
$4$ & $31^3$ & $2.81$e$+00$ ($40$\%) & $3.37$e$-01$ ($78$\%) & $10$ ($10$) & $3.67$e$+00$ ($79$\%) \\
$8$ & $63^3$ & $1.64$e$+01$ ($28$\%) & $2.24$e$+00$ ($68$\%) & $10$ ($10$) & $2.52$e$+01$ ($70$\%) \\
$16$ & $127^3$ & $1.12$e$+02$ ($19$\%) & $1.64$e$+01$ ($60$\%) & $9$ ($9$) & $1.82$e$+02$ ($65$\%) \\
$32$ & $255^3$ & $8.22$e$+02$ ($11$\%) & $1.24$e$+02$ ($50$\%) & $9$ ($9$) & $1.39$e$+03$ ($52$\%) \\
\hline
\hline
\end{tabular}

\caption[Numerical results for velocity field \eqref{LS:vel:3D_4} in 3D with the recursive approach.]{Numerical results for velocity field \eqref{LS:vel:3D_4} in 3D with the recursive approach.}
\label{LS:tab:3D_recursive_4}
\end{table}

From the numerical tests we see that, same as the 2D cases, the iteration numbers remain
essentially independent of the problem size. The preconditioner converges in a few iterations for
all the test cases. Another highlight is that, the recursive approach requires only zero or one more
iteration compared to the non-recursive approach, which means that the recursive sweeping
factorization for the quasi-2D linear systems keeps the total approximation error almost at the same level.

\section{Sparsifying scheme as a direct method}
\label{sec:stencil}

In this section, we show that the compact stencils acquired by the sparsifying scheme can be viewed as accurate discretizations of the Helmholtz equation. Specifically, we will solve the 2D homogeneous Helmholtz equation with the compact scheme introduced in the sparsification process, and compare it with the Quasi-Stabilized FEM (QSFEM) method in \cite{babuska1995}. As we shall see from the numerical tests, both methods did comparably well at minimizing the pollution error with only a small number of points per wavelength.

\begin{figure}
[!ht]
\centering
\includegraphics[width=0.38\textwidth]{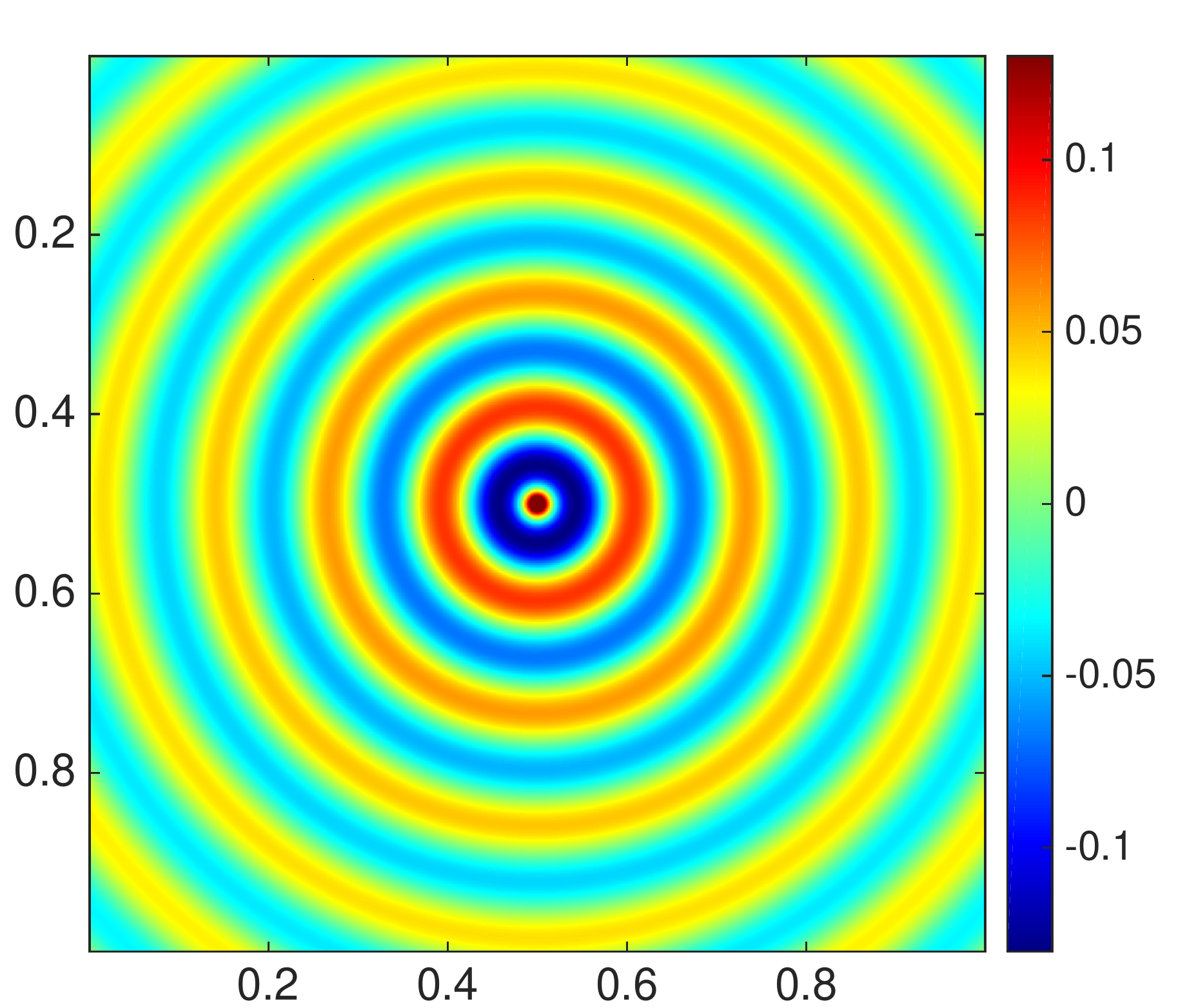}
\caption{An example of the 2D Green's function centered at $(0.5, 0.5)$ with $\omega/(2\pi)=8$.}
\label{fig:2D_Green}
\end{figure}

Let's consider
\begin{equation}
\label{eqn:2D_Helm_comarison}
(-\Delta - \omega^2) u(x) = f(x), \quad x \in \bR^2,
\end{equation}
where $f(x)$ is a delta source centered at $(0.5, 0.5)$. The exact solution is given by the Green's function with a shift of the center (see Figure \ref{fig:2D_Green} for an example).

For the sparsifying scheme, we have the discrete equation
\begin{gather*}
\alpha^* u_{\mu_i} = \beta^* f_{\mu_i}
\end{gather*}
for each of the interior point $i$, where $\alpha$ and $\beta$ are 9-point stencils given by \eqref{eqn:alpha} and \eqref{eqn:beta} respectively, and $f$ is the discrete delta function.

For the QSFEM method, the 9-point stencil for $u$ is given by
\begin{gather*}
A
=
\begin{bmatrix}
A_2 & A_1 & A_2 \\
A_1 & A_0 & A_1 \\
A_2 & A_1 & A_2
\end{bmatrix},
\end{gather*}
where
\begin{align*}
A_0 & = 4,\\
A_1 & = 2 \dfrac{c_1(\kappa)s_1(\kappa) - c_2(\kappa) s2(\kappa)}{c_2(\kappa) s_2(\kappa)(c_1(\kappa)+s_1(\kappa)) - c_1(\kappa) s_1(\kappa) (c_2(\kappa) + s_2(\kappa))},\\
A_2 & = 2 \dfrac{c_2(\kappa) + s_2(\kappa) - c_1(\kappa) - s1(\kappa)}{c_2(\kappa) s_2(\kappa)(c_1(\kappa)+s_1(\kappa)) - c_1(\kappa) s_1(\kappa) (c_2(\kappa) + s_2(\kappa))},\\
c_1(\kappa) &\coloneqq \cos\left(\kappa \cos \dfrac{\pi}{16}\right),
\quad
s_1(\kappa) \coloneqq \cos\left(\kappa \sin \dfrac{\pi}{16}\right), \\
c_2(\kappa) &\coloneqq \cos\left(\kappa \cos \dfrac{3\pi}{16}\right),
\quad
s_2(\kappa) \coloneqq \cos\left(\kappa \sin \dfrac{3\pi}{16}\right), \\
\kappa & \coloneqq \omega h,
\end{align*}
and $h$ is the step size. The right-hand side is the discrete delta function with a scaling.

We solve the 2D homogeneous Helmholtz equation \eqref{eqn:2D_Helm_comarison} and compare the phase errors against the true solution. Specifically, we write the solutions $u(x)$ as $u(x) = A(x) e^{2\pi \ii \phi(x)}$ and compare the phase $\phi(x)$ with the one acquired by the Green's function. The boundary points are discretized by a slowly turning-up PML such that the reflection error is negligible compared to the phase error.

\begin{figure}
[!ht]
\vspace{-25pt}
\centering
\subfigure[Sparsifying scheme at 5 p.p.w.]{
\includegraphics[width=0.38\textwidth]{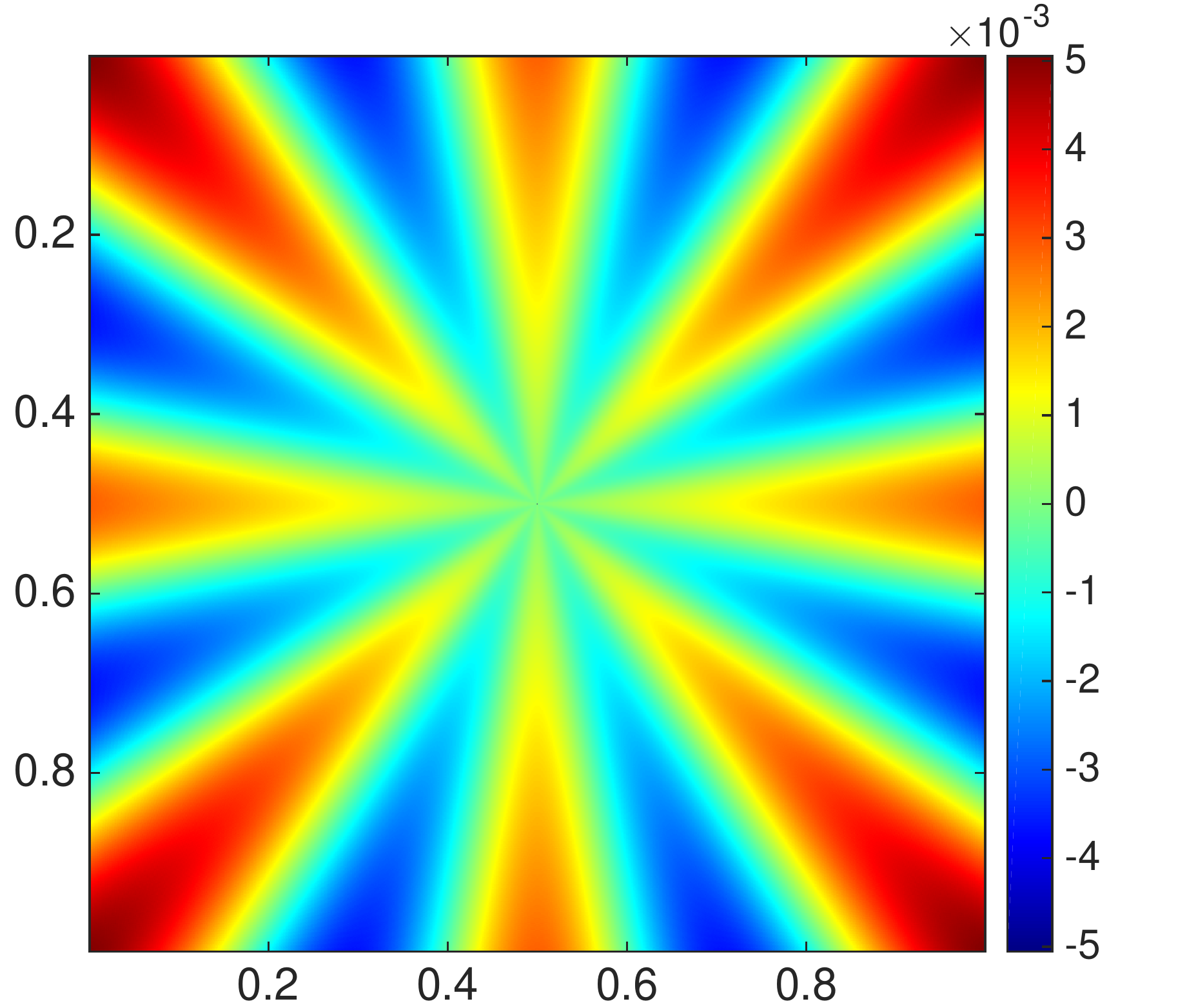}
}
\subfigure[QSFEM at 5 p.p.w.]{
\includegraphics[width=0.38\textwidth]{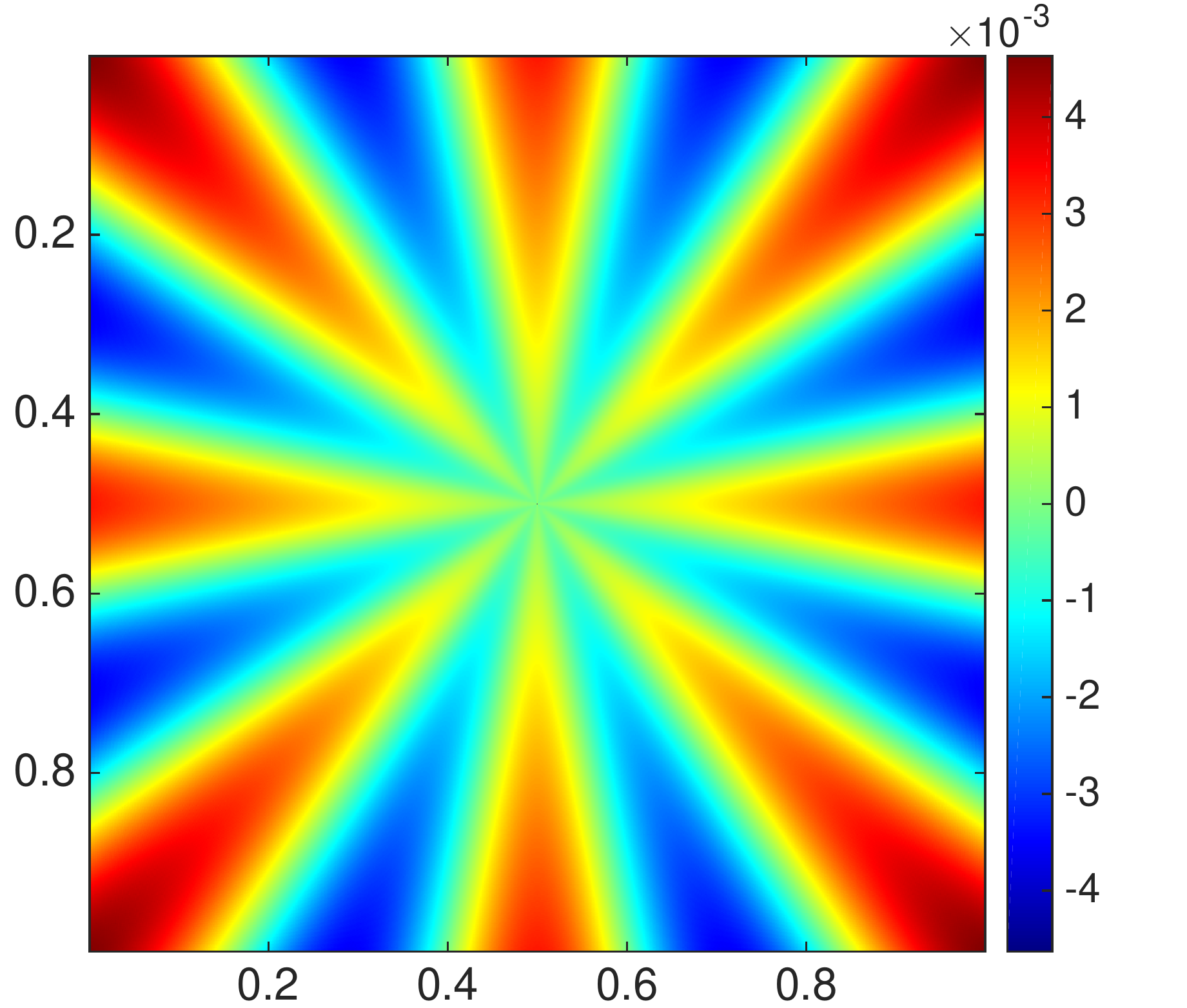}
}

\vspace{-5pt}
\subfigure[Sparsifying scheme at 4 p.p.w.]{
\includegraphics[width=0.38\textwidth]{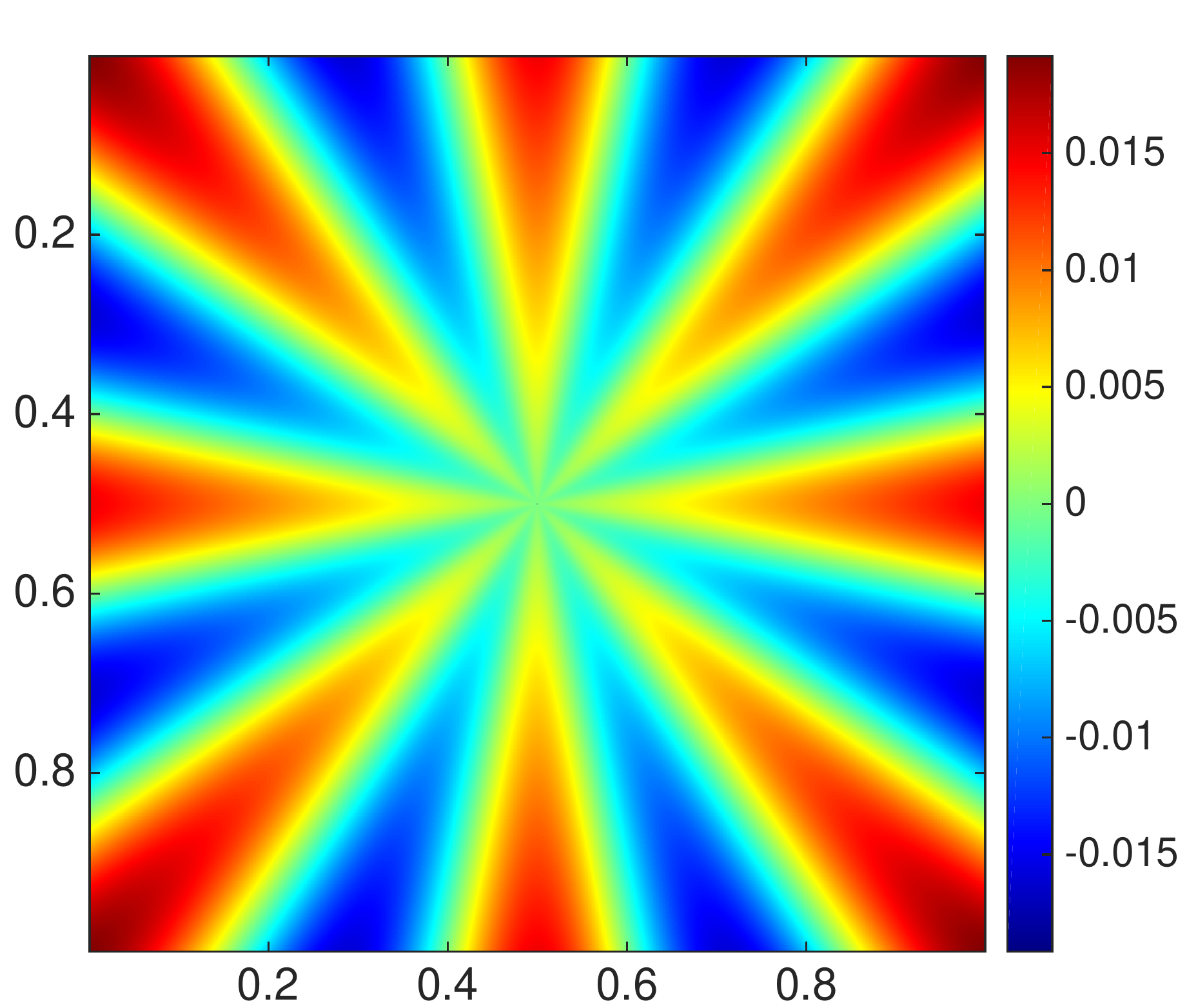}
}
\subfigure[QSFEM at 4 p.p.w.]{
\includegraphics[width=0.38\textwidth]{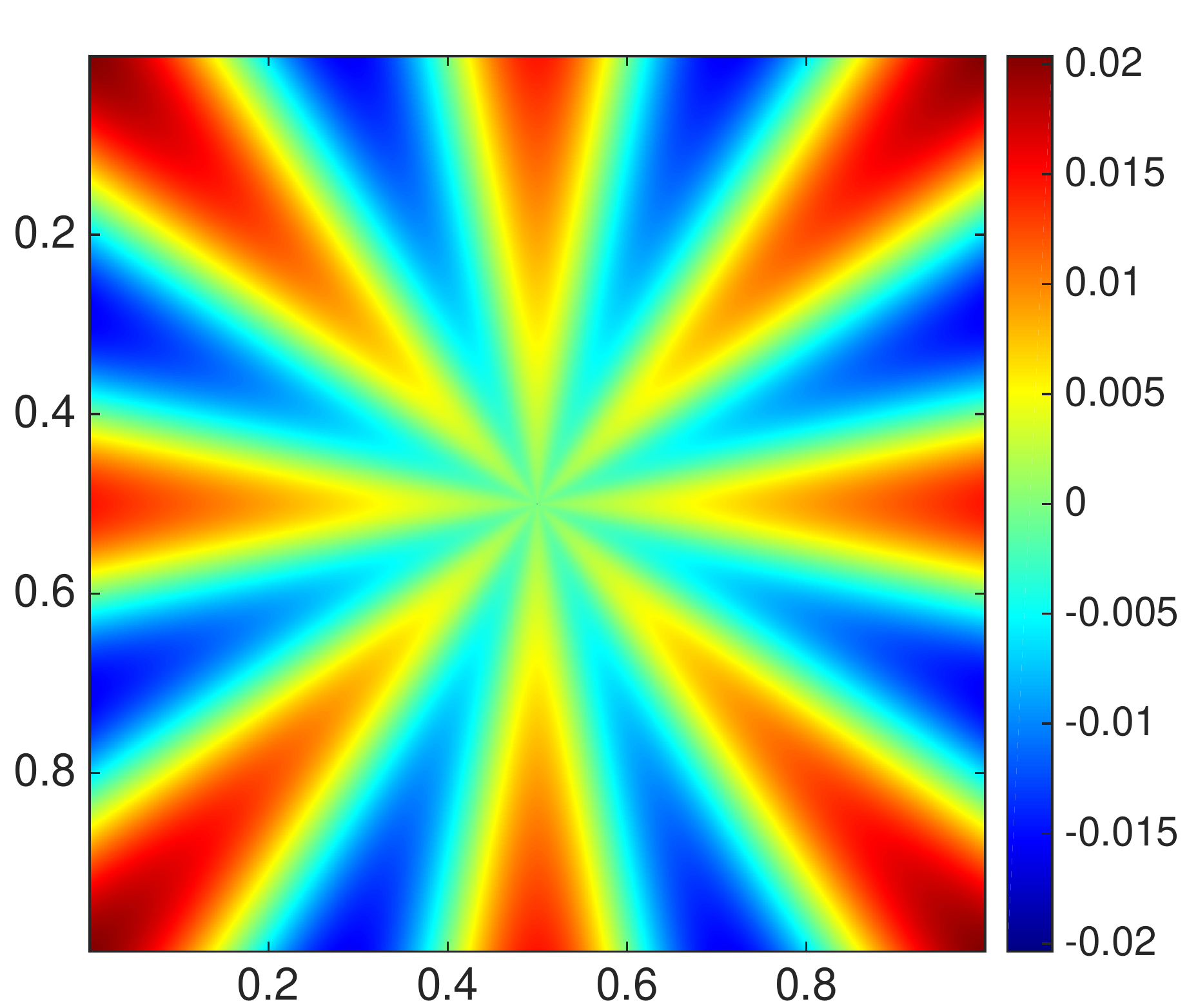}
}

\vspace{-5pt}
\subfigure[Sparsifying scheme at 3 p.p.w.]{
\includegraphics[width=0.38\textwidth]{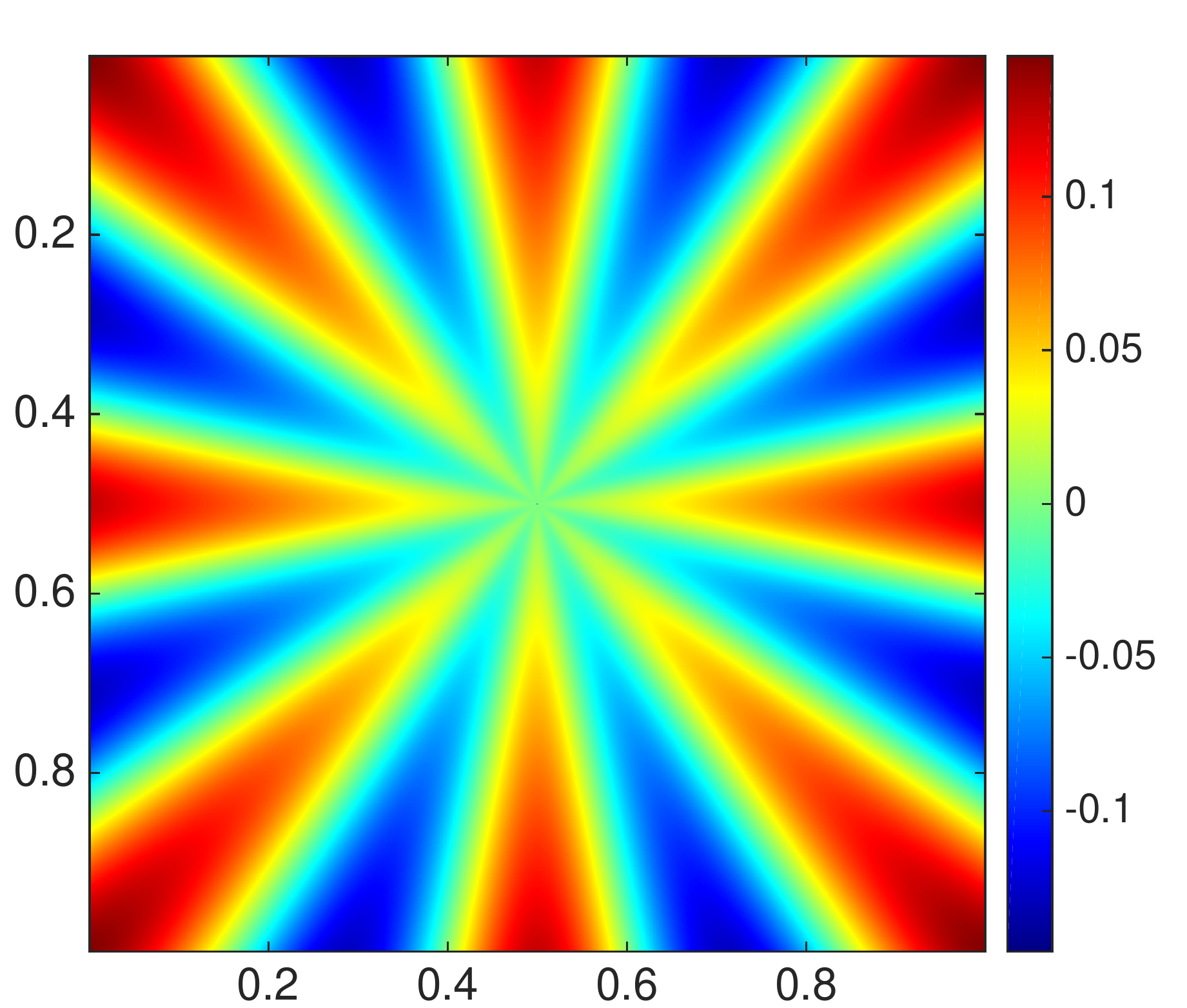}
}
\subfigure[QSFEM at 3 p.p.w.]{
\includegraphics[width=0.38\textwidth]{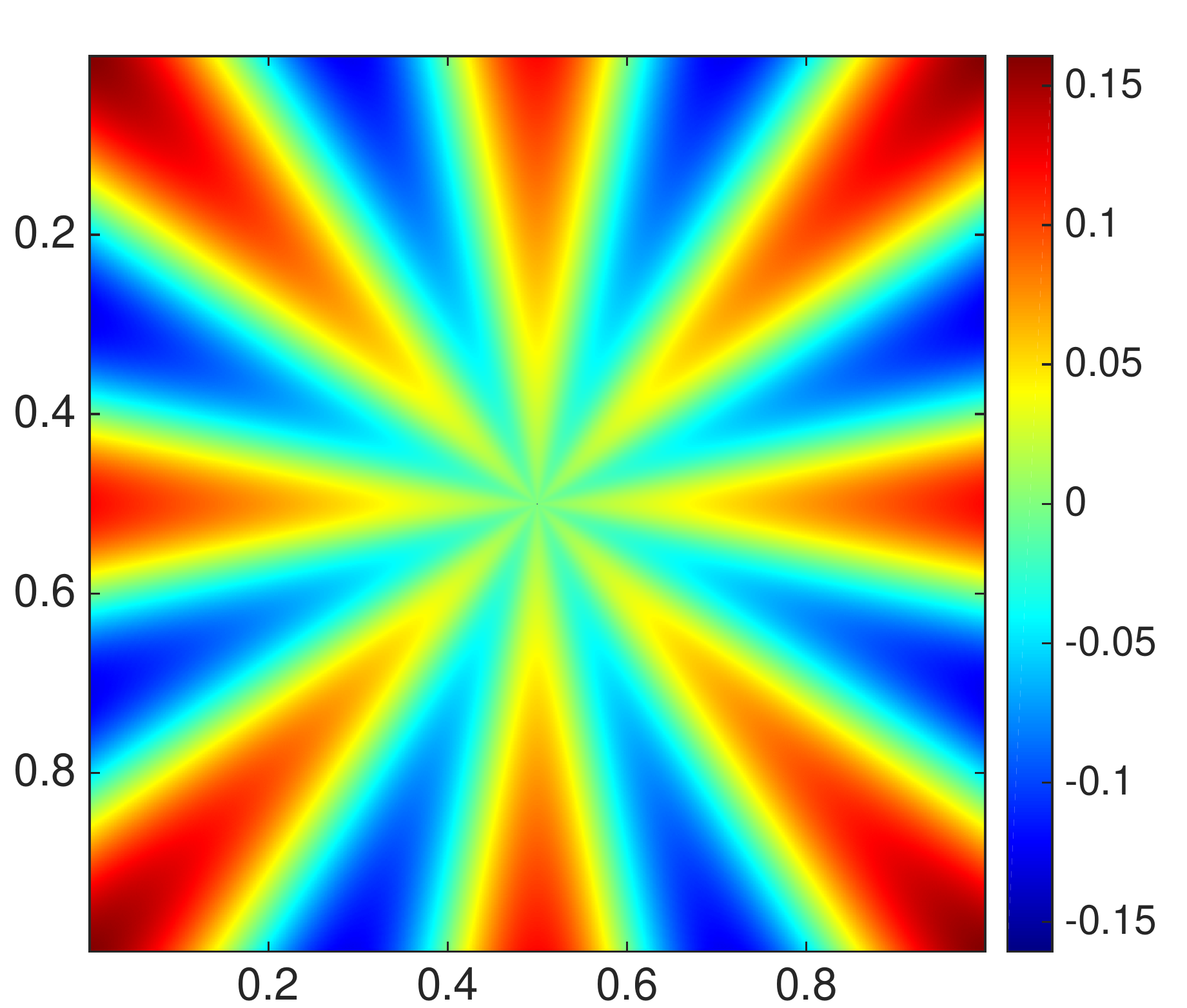}
}
\caption{The phase errors of the two schemes for $\omega/(2\pi)=1024$. Phase error is defined as $\phi(x) - \phi_{G}(x)$ where $\phi(x)$ is computed from the numerical solution $u(x) = A(x) e^{2\pi \ii \phi(x)}$, and $\phi_{G}(x)$ is acquired by the phase of the Green's function (exact solution). From top to bottom are test cases for 5, 4 and 3 points per wavelength respectively. We see that both methods behave similarly in terms of the phase error. For the hardest test cases (3 p.p.w with 1024 waves across each dimension), the phase shifts at the far field (four corners) are about 1/6. This corresponds to about $2.3\times 10^{-4}$ relative phase error. In other words, the distortion of each wavelength of both methods are on the order of $10^{-4}$ with only three points per wavelength.}
\label{fig:2D_comparison}
\end{figure}

Figure \ref{fig:2D_comparison} shows the phase errors for a large test case (1024 waves across each dimension) with a small number of points (3 to 5) per wavelength. From the tests we see that the phase error of the sparsifying scheme is comparable to the one of the QSFEM method in \cite{babuska1995}.

We would like to comment that, as pointed out in \cite{Babuska1997}, 2D compact stencils can be optimized to reduce the pollution error, but cannot completely eliminate it. For example, in Figure \ref{fig:2D_comparison}, the phase shifts near the four corners are about 1/6 for the 3 p.p.w. test cases, which is not negligible for practical usage. Hence for large problems, one would eventually have to increase the stencil width, or use more points per wavelength.

\section{Conclusions and future work}
\label{sec:conclusion}

This paper presents the sparsify-and-sweep preconditioner for the Lippmann-Schwinger equation in 2D
and 3D. The preconditioner involves two steps. The first step is to sparsify the system by
introducing the compact stencil sparsifying scheme. The second step is to apply the sweeping
factorization to the sparsified system. Numerical results show that the iteration number is
essentially independent of the angular frequency $\omega$.

Though the cost is reduced to linear, potential improvements can be made regarding
parallelizations. First, the factorization of the auxiliary subproblems are completely independent,
thus can be done in parallel, especially when the recursive approach in 3D is adopted where there
are $O(n^2)$ quasi-1D subproblems that can be processed at the same time. Second, the setup and
application processes of the nested dissection algorithm can also be parallelized for independent
skeleton fronts (see \cite{poulson2013parallel} for example). Third, the two sweeping fronts during
the application process are also independent and can be processed in a parallel way.

Another future work is on the sparsification of dense systems by the data-fitting approach. This
approach was first proposed by Ying in \cite{Ying2015, Ying2015sp} for solving highly indefinite
systems including time-independent high frequency wave propagations with radiation conditions or
periodic boundary conditions. This paper generalizes it to incorporate the PML approach. There have been some explorations and applications of this sparsification method,
such as solving the nonlinear eigenvalue problems in soliton systems \cite{Lu2016soliton}. This
data-fitting approach to design local schemes is quite different from most classical approaches, and
could be potentially generalized to other types of integral equations and dense systems.

\section*{Acknowledgments}
The authors are partially supported by the National Science
Foundation under award DMS-1521830 and the U.S. Department of Energy's
Advanced Scientific Computing Research program under award
DE-FC02-13ER26134/DE-SC0009409.

\bibliographystyle{abbrv}
\bibliography{references}

\begin{thebibliography}{10}

\bibitem{greengard2016}
S.~Ambikasaran, C.~Borges, L.-M. Imbert-Gerard, and L.~Greengard.
\newblock {Fast, Adaptive, High-Order Accurate Discretization of the
  Lippmann--Schwinger Equation in Two Dimensions}.
\newblock {\em SIAM Journal on Scientific Computing}, 38(3):A1770--A1787, 2016.

\bibitem{Andersson2005}
F.~Andersson and A.~Holst.
\newblock {A Fast, Bandlimited Solver for Scattering Problems in Inhomogeneous
  Media}.
\newblock {\em Journal of Fourier Analysis and Applications}, 11(4):471--487,
  2005.

\bibitem{asvadurov2003optimal}
S.~Asvadurov, V.~Druskin, M.~N. Guddati, and L.~Knizhnerman.
\newblock {On optimal finite-difference approximation of PML}.
\newblock {\em SIAM Journal on Numerical Analysis}, 41(1):287--305, 2003.

\bibitem{Babuska1997}
I.~M. Babu{\v{s}}ka and S.~A. Sauter.
\newblock {Is the Pollution Effect of the FEM Avoidable for the Helmholtz
  Equation Considering High Wave Numbers?}
\newblock {\em SIAM Journal on Numerical Analysis}, 34(6):2392--2423, 1997.

\bibitem{babuska1995}
I.~Babuška, F.~Ihlenburg, E.~T. Paik, and S.~A. Sauter.
\newblock A generalized finite element method for solving the helmholtz
  equation in two dimensions with minimal pollution.
\newblock {\em Computer Methods in Applied Mechanics and Engineering},
  128(3):325 -- 359, 1995.

\bibitem{Berenger1994}
J.-P. Berenger.
\newblock {A perfectly matched layer for the absorption of electromagnetic
  waves}.
\newblock {\em J. Comput. Phys.}, 114(2):185--200, 1994.

\bibitem{Bruno2004}
O.~P. Bruno and E.~Hyde.
\newblock {An efficient, preconditioned, high-order solver for scattering by
  two-dimensional inhomogeneous media}.
\newblock {\em Journal of Computational Physics}, 200(2):670 -- 694, 2004.

\bibitem{Chen2002}
Y.~Chen.
\newblock {A Fast, Direct Algorithm for the Lippmann-Schwinger Integral
  Equation in Two Dimensions}.
\newblock {\em Advances in Computational Mathematics}, 16(2):175--190, 2002.

\bibitem{Chen2013a}
Z.~Chen and X.~Xiang.
\newblock {A source transfer domain decomposition method for Helmholtz
  equations in unbounded domain}.
\newblock {\em SIAM J. Numer. Anal.}, 51(4):2331--2356, 2013.

\bibitem{Chen2013b}
Z.~Chen and X.~Xiang.
\newblock {A source transfer domain decomposition method for Helmholtz
  equations in unbounded domain Part II: Extensions}.
\newblock {\em Numer. Math. Theory Methods Appl.}, 6(3):538--555, 2013.

\bibitem{Chew1994}
W.~C. Chew and W.~H. Weedon.
\newblock {A 3D perfectly matched medium from modified Maxwell's equations with
  stretched coordinates}.
\newblock {\em Microw. Opt. Techn. Let.}, 7(13):599--604, 1994.

\bibitem{Duan2009}
R.~Duan and V.~Rokhlin.
\newblock {High-order Quadratures for the Solution of Scattering Problems in
  Two Dimensions}.
\newblock {\em J. Comput. Phys.}, 228(6):2152--2174, Apr. 2009.

\bibitem{Ying2011a}
B.~Engquist and L.~Ying.
\newblock {Sweeping preconditioner for the Helmholtz equation: hierarchical
  matrix representation}.
\newblock {\em Comm. Pure Appl. Math.}, 64(5):697--735, 2011.

\bibitem{Ying2011b}
B.~Engquist and L.~Ying.
\newblock {Sweeping preconditioner for the Helmholtz equation: moving perfectly
  matched layers}.
\newblock {\em Multiscale Model. Simul.}, 9(2):686--710, 2011.

\bibitem{george1973nested}
A.~George.
\newblock {Nested dissection of a regular finite element mesh}.
\newblock {\em SIAM J. Numer. Anal.}, 10:345--363, 1973.
\newblock Collection of articles dedicated to the memory of George E. Forsythe.

\bibitem{Johnson2008}
S.~G. Johnson.
\newblock {Notes on Perfectly Matched Layers (PMLs)}.
\newblock {\em Lecture notes, Massachusetts Institute of Technology,
  Massachusetts}, 2008.

\bibitem{kristek2009brief}
J.~Kristek, P.~Moczo, and M.~Galis.
\newblock {A brief summary of some PML formulations and discretizations for the
  velocity-stress equation of seismic motion}.
\newblock {\em Studia Geophysica et Geodaetica}, 53(4):459, 2009.

\bibitem{Lanzara2004}
F.~Lanzara, V.~Maz'ya, and G.~Schmidt.
\newblock {Numerical Solution of the Lippmann-Schwinger equation by approximate
  approximations}.
\newblock {\em Journal of Fourier Analysis and Applications}, 10(6):645--660,
  2004.

\bibitem{Liu2016}
F.~Liu and L.~Ying.
\newblock {Additive sweeping preconditioner for the Helmholtz equation}.
\newblock {\em Multiscale Modeling \& Simulation}, 14(2):799--822, 2016.

\bibitem{Liu2015}
F.~Liu and L.~Ying.
\newblock {Recursive Sweeping Preconditioner for the Three-Dimensional
  Helmholtz Equation}.
\newblock {\em SIAM Journal on Scientific Computing}, 38(2):A814--A832, 2016.

\bibitem{liu2016localized}
F.~Liu and L.~Ying.
\newblock {Localized sparsifying preconditioner for periodic indefinite
  systems}.
\newblock {\em Commun. Math. Sci.}, 15(4):1155--1169, 2017.

\bibitem{Lu2016soliton}
J.~Lu and L.~Ying.
\newblock {Sparsifying preconditioner for soliton calculations}.
\newblock {\em Journal of Computational Physics}, 315:458 -- 466, 2016.

\bibitem{poulson2013parallel}
J.~Poulson, B.~Engquist, S.~Li, and L.~Ying.
\newblock {A parallel sweeping preconditioner for heterogeneous 3D Helmholtz
  equations}.
\newblock {\em SIAM Journal on Scientific Computing}, 35(3):C194--C212, 2013.

\bibitem{Sifuentes2010}
J.~Sifuentes.
\newblock {\em {Preconditioned iterative methods for inhomogeneous acoustic
  scattering applications}}.
\newblock PhD thesis, Rice University Houston, 2010.

\bibitem{Stolk2013}
C.~C. Stolk.
\newblock {A rapidly converging domain decomposition method for the Helmholtz
  equation}.
\newblock {\em J. Comput. Phys.}, 241(0):240 -- 252, 2013.

\bibitem{stolk2016}
C.~C. Stolk.
\newblock A dispersion minimizing scheme for the 3-d helmholtz equation based
  on ray theory.
\newblock {\em Journal of Computational Physics}, 314(Supplement C):618 -- 646,
  2016.

\bibitem{Vainikko2000}
G.~Vainikko.
\newblock {Fast solvers of the Lippmann-Schwinger equation}.
\newblock In {\em Direct and inverse problems of mathematical physics}, pages
  423--440. Springer, 2000.

\bibitem{vico2016}
F.~Vico, L.~Greengard, and M.~Ferrando.
\newblock Fast convolution with free-space green's functions.
\newblock {\em Journal of Computational Physics}, 323(Supplement C):191 -- 203,
  2016.

\bibitem{Vion2014}
A.~Vion and C.~Geuzaine.
\newblock {Double sweep preconditioner for optimized Schwarz methods applied to
  the Helmholtz problem}.
\newblock {\em J. Comput. Phys.}, 266(0):171 -- 190, 2014.

\bibitem{Ying2015sp}
L.~Ying.
\newblock {Sparsifying Preconditioner for Pseudospectral Approximations of
  Indefinite Systems on Periodic Structures}.
\newblock {\em Multiscale Modeling \& Simulation}, 13(2):459--471, 2015.

\bibitem{Ying2015}
L.~Ying.
\newblock {Sparsifying Preconditioner for the Lippmann-Schwinger Equation}.
\newblock {\em Multiscale Modeling \& Simulation}, 13(2):644--660, 2015.

\bibitem{Demanet2014}
L.~Zepeda-N{\'u}{\~n}ez and L.~Demanet.
\newblock {The method of polarized traces for the 2D Helmholtz equation}.
\newblock {\em Journal of Computational Physics}, 308:347 -- 388, 2016.

\bibitem{Leonardo2016fast}
L.~Zepeda-N{\'u}{\~n}ez and H.~Zhao.
\newblock {Fast Alternating BiDirectional Preconditioner for the 2D
  High-Frequency Lippmann-Schwinger Equation}.
\newblock {\em SIAM Journal on Scientific Computing}, 38(5):B866--B888, 2016.

\end{thebibliography}
\end{document}